\documentclass[12pt,leqno,twoside]{article}
\usepackage[latin1]{inputenc}
\usepackage[english]{babel}
\usepackage{amsthm,amssymb,amsmath,graphics,enumerate,epsfig,fancyhdr,url}

\newcommand{\ThisTitle}{}
\title{\ThisTitle}
\author{Magnus Dehli Vigeland\thanks{Department of Mathematics, University of Oslo, Norway.
{\it Email}\,: {\tt magnusv@math.uio.no}}
}
\date{}

\renewcommand{\ThisTitle}{Smooth tropical surfaces with infinitely many tropical lines}

\theoremstyle{plain}
\newtheorem{theo}{Theorem}[section]
\newtheorem{prop}[theo]{Proposition}
\newtheorem{lem}[theo]{Lemma}
\newtheorem{cor}[theo]{Corollary}
\newtheorem*{theor}{Theorem}

\theoremstyle{definition}
\newtheorem{dfn}[theo]{Definition}

\theoremstyle{remark}
\newtheorem{ex}[theo]{Example}
\newtheorem{remark}[theo]{Remark}

\newcommand{\pil}{\rightarrow}
\newcommand{\nn}{\mathbb{N}}
\newcommand{\zz}{\mathbb{Z}}

\newcommand{\rr}{\mathbb{R}}
\newcommand{\La}{\Lambda}
\newcommand{\la}{\lambda}

\newcommand{\Ga}{\Gamma}

\newcommand{\dd}{\mathcal{D}}
\newcommand{\cc}{\mathcal{C}}
\newcommand{\ks}{\mathcal{S}}
\newcommand{\kt}{\mathcal{T}}
\newcommand{\kg}{\mathcal{G}}
\newcommand{\ke}{\mathcal{E}}
\newcommand{\ff}{\mathcal{F}}
\newcommand{\vv}{\mathcal{V}}
\newcommand{\qq}{\mathcal{Q}}
\newcommand{\skrifta}{\mathcal{A}}

\newcommand{\sub}{\subseteq}
\newcommand{\wb}[1]{\overline{#1}}
\newcommand{\wvec}[1]{\overrightarrow{#1}}
\newcommand{\inner}[2]{\langle #1,#2\rangle}

\DeclareMathOperator{\conv}{conv}
\DeclareMathOperator{\inter}{int}

\DeclareMathOperator{\vol}{vol}
\DeclareMathOperator{\Aff}{Aff}
\DeclareMathOperator{\Subdiv}{Subdiv}
\DeclareMathOperator{\im}{im}
\DeclareMathOperator{\Fac}{Fac}

\begin{document}
\setcounter{page}{1}

\maketitle
\begin{abstract}
We study the tropical lines contained in smooth tropical surfaces in $\rr^3$. On smooth tropical quadric surfaces we find two one-dimensional families of tropical lines, like in classical algebraic geometry. Unlike the classical case, however, there exist smooth tropical surfaces of any degree with infinitely many tropical lines. 
\end{abstract}

\section{Introduction}
Tropical geometry has during the last few years become an increasingly popular field of mathematics. This is not least due to the many fascinating similarities with classical geometry. In this paper we examine tropical analogues of the following well-known results in classical algebraic geometry:

\begin{enumerate}[(I)]
\item Any smooth quadric surface has two rulings of lines,
\item Any smooth surface of degree greater than two, has at most finitely many lines.
\end{enumerate}

While a lot of work has been done lately on tropical curves (e.g. \cite{Mikh, mikh06, mikh07, GM2, jinv, vig04}), comparatively little is known in higher dimensions. 
The usual way of defining a tropical variety is as the {\em tropicalization} of an algebraic variety defined over an algebraically closed field with a non-Archimedean valuation (see e.g. \cite{RGST}). 
In the case of hypersurfaces, however, a more inviting, geometric definition is possible. For example, a tropical surface in $\rr^3$ is precisely the non-linear locus of a continuous convex piecewise linear function \mbox{$f:\rr^3\pil\rr$} with rational slopes. It is an unbounded two-dimensional polyhedral complex, with {\em zero tension} at each $1$-cell. Furthermore, it is dual to a regular subdivision of the Newton polytope of $f$ (when $f$ is regarded as a {\em tropical polynomial}). The tropical surface is {\em smooth} if this subdivision is an elementary (unimodular) triangulation.

Tropical varieties of higher codimension are in general more difficult to grasp. However, the only such varieties we are interested in here, are tropical lines in $\rr^3$. These were given an explicit geometric description in \cite{RGST}, on which we base our definition. 
As an analogue of (I) above, we prove that:
\begin{theor}
Any smooth tropical quadric surface $X$ has a unique compact 2-cell $\wb{X}$. For any point $p\in\wb{X}$, there exist two tropical lines on $X$ containing $p$.
\end{theor}

While in classical geometry, any two distinct points in $\rr^3$ lie on a unique line, this is only true generically for tropical lines. In fact, for special choices of $p,q\in\rr^3$ there are infinitely many tropical lines containing $p$ and $q$. We show that such families of tropical lines can also exist on a smooth tropical surface. As a consequence, we get the following result, in contrast to (II) above:

\begin{theor}
There exist tropical surfaces of any degree, with infinitely many tropical lines.
\end{theor}

The paper is organized as follows: In sections \ref{section:pol} and \ref{section:trop} we give some necessary background on convex geometry and tropical geometry, respectively. In particular, the concept of a {\em two-point family} of tropical lines in $\rr^3$ is defined in Section \ref{section:troplines}. Then follows two technical sections, \ref{section:construction} and \ref{section:exits}. The former of these deals with constructions of regular elementary triangulations, while the latter contains an analysis of certain lattice polytopes. In Section \ref{section:genprop} we explore the general properties of tropical lines contained in smooth tropical surfaces, and in Section \ref{section:quad} we use these to study tropical lines on quadric surfaces. Section \ref{twopoint} concerns two-point families of tropical lines on smooth tropical surfaces. Finally, Section \ref{section:higher} contains our results for tropical surfaces of higher degrees.

\section{Lattice polytopes and subdivisions}\label{section:pol}
\subsection{Convex polyhedra and polytopes}
A {\em convex polyhedron} in $\rr^n$ is the intersection of finitely many closed halfspaces. A {\em cone} is a convex polyhedron, all of whose defining hyperplanes contain the origin. A {\em convex polytope} is a bounded convex polyhedron. Equivalently, a convex polytope can be defined as the convex hull of a finite set of points in $\rr^n$. Throughout this paper, all polyhedra and polytopes will be assumed to be convex unless explicitly stated otherwise.

For any polyhedron $\Delta\sub\rr^n$ we denote its affine hull by $\Aff(\Delta)$, and its relative interior (as a subset of $\Aff(\Delta)$) by $\inter(\Delta)$. The dimension of $\Delta$ is defined as $\dim\Aff(\Delta)$. By convention, $\dim \emptyset=-1$. A {\em face} of $\Delta$ is a polyhedron of the form $\Delta\cap H$, where $H$ is a hyperplane such that $\Delta$ is entirely contained in one of the closed halfspaces defined by $H$. In particular, the empty set is considered a face of $\Delta$. Faces of dimensions $0$, $1$ and $n-1$ are called {\em vertices}, {\em edges} and {\em facets} of $\Delta$, respectively. If $\Delta$ is a polytope, then the vertices of $\Delta$ form the minimal set $\skrifta$ such that $\Delta=\conv(\skrifta)$. 

Let $F$ be a facet of a polyhedron $\Delta\sub\rr^n$, where $\dim \Delta\leq n$. A vector $v$ is {\em pointing inwards} (resp. {\em pointing outwards})
 from $F$ relative to $\Delta$ if, for some positive constant $t$, the vector $tv$ (resp. $-tv$) starts in $F$ and ends in $\Delta\smallsetminus F$. If in addition $v$ is orthogonal to $F$, $v$ is an {\em inward normal vector} (resp. {\em outward normal vector}) of $F$ relative to $\Delta$. Using the notation $\inner{\:}{\:}$ for the Euclidean inner product, a straightforward consequence of these definitions is:
\begin{lem}\label{inout}
Let $\Delta$ and $F$ be as above, and let $v$ be a vector. The following are equivalent:
\begin{enumerate}
\item $v$ is an inward (resp. outward) normal vector of $F$ relative to $\Delta$,
\item $v$ is orthogonal to $F$, and $\inner{u}{v}>0$ (resp. $<0$) for some vector $u$ pointing inwards from $F$ relative to $\Delta$,
\item $\inner{u}{v}>0$ (resp. $<0$) for all vectors $u$ pointing inwards from $F$ relative to $\Delta$.
\end{enumerate}
\end{lem}


If all the vertices of $\Delta$ are contained in $\zz^n$, we call $\Delta$ a {\em lattice polyhedron}, or {\em lattice polytope} if it is bounded. A lattice polytope in $\rr^n$ is {\em primitive} if it contains no lattice points other than its vertices. It is {\em elementary} (or {\em unimodular}) if it is $n$-dimensional and its volume is $\frac1{n!}$. 
Obviously, every elementary polytope is also primitive, while the other implication is not true in general. For instance, the unit square in $\rr^2$ is primitive, but not elementary.
 
Most of the polytopes we are interested in will be simplices, i.e., the convex hull of $n+1$ affinely independent points. In $\rr^2$, the primitive simplices are precisely the elementary ones, namely the lattice triangles of area $\frac12$. (This is an immediate consequence of Pick's theorem.) In higher dimensions, the situation is very different: There is no upper limit for the volume of a primitive simplex in $\rr^n$, when $n\geq 3$. The standard example of this is the following: Let $p,q\in \nn$ be relatively prime, with $p<q$, and let $T_{p,q}$ be the tetrahedron with vertices in $(0,0,0)$, $(1,0,0)$, $(0,1,0)$ and $(1,p,q)$. Then $T_{p,q}$ is a primitive simplex of volume $\frac{q}6$.

\subsection{Polyhedral complexes and subdivisions}\label{subsec:subs}
A (finite) {\em polyhedral complex} in $\rr^n$ is a finite collection $X$ of convex polyhedra, called {\em cells}, such that
\begin{itemize}
\item if $C\in X$, then all faces of $C$ are in $X$, and
\item if $C,C'\in X$, then $C\cap C'$ is a face of both $C$ and $C'$.
\end{itemize}
The $d$-dimensional elements of $X$ are called the {\em $d$-cells} of $X$. The dimension of $X$ itself is defined as $\max\{\dim C\:|\: C\in X\}$. Furthermore, if all the maximal cells (w.r.t. inclusion) have the same dimension, we say that $X$ is of {\em pure dimension}.

A polyhedral complex, all of whose cells are cones, is a {\em fan}.

A {\em subdivision} of a polytope $\Delta$ is a polyhedral complex $\ks$ such that $|\ks|=\Delta$, where $|\ks|$ denotes the union of all the elements of $\ks$. It follows that $\ks$ is of pure dimension $\dim \Delta$. If all the maximal elements of $\ks$ are simplices, we call $\ks$ a {\em triangulation}. If $\ks$ and $\ks'$ are subdivisions of the same polytope, we say that $\ks'$ is a {\em refinement} of $\ks$ if for all $C'\in \ks'$ there is a $C\in\ks$ such that
$C'\sub C$.

If $\Delta$ is a lattice polytope, we can consider {\em lattice subdivisions} of $\Delta$, i.e., subdivisions in which every element is a lattice polytope. In particular, a lattice subdivision is {\em primitive} (resp. {\em elementary}) if all its maximal elements are primitive (resp. elementary). We write down some noteworthy properties of these subdivisions:
\begin{itemize}
\item Every elementary subdivision is a primitive triangulation.
\item In a primitive subdivision, all elements (not only the maximal) are primitive.
\item For any lattice polytope, its lattice subdivisions with no non-trivial refinements are precisely its primitive triangulations.
\end{itemize}
\subsection{Regular subdivisions and the secondary fan}\label{section:reg}
Let $\Delta=\conv(\skrifta)$ where $\skrifta$ is a finite set of points in $\rr^n$. Any function $\alpha\colon \skrifta\pil \rr$ will induce a lattice subdivision of $\Delta$ in the following way. Consider the polytope
\begin{equation*}
\conv(\{(v,\alpha(v))\:|\: v\in\skrifta\})\in\rr^{n+1}.
\end{equation*}
Projecting the top faces of this polytope to $\rr^n$, forgetting the last coordinate, gives a collection of subpolytopes of $\Delta$. They form a subdivision $\ks_\alpha$ of $\Delta$. The function $\alpha$ is called a {\em lifting function} associated to $\ks_\alpha$. 

\begin{dfn} 
A lattice subdivision $\ks$ of $\conv(\skrifta)$ is {\em regular} if $\ks=\ks_\alpha$ for some $\alpha\colon\skrifta\pil\rr$.
\end{dfn}

The set of regular subdivisions of $\conv(\skrifta)$ has an interesting geometric structure, as observed by Gelfand, Kapranov and Zelevinsky in \cite{GKZ}. 
Suppose $\skrifta\sub\rr^n$ consists of $k$ points. For a fixed ordering of the points in $\skrifta$, the space $\rr^\skrifta\simeq \rr^k$ is a parameter space for all functions $\alpha\colon\skrifta\pil\rr$. For a given given regular subdivision $\ks$ of $\conv(\skrifta)$, let $K(\ks)$ be the set of all functions $\alpha\in \rr^\skrifta$ which induce $\ks$. The following is proved in \cite[Chapter 7]{GKZ}:
\begin{prop}\label{prop:secfan}
Let $\ks$ and $\ks'$ be any regular subdivisions of $\conv(\skrifta)$. Then:
\begin{enumerate}
\item $K(\ks)$ is a cone in $\rr^\skrifta$.
\item $\ks'$ is a refinement of $\ks$ if and only if $K(\ks)$ is a face of $K(\ks')$.
\item The cones $\{K(\ks)\:|\:\text{$\ks$ is a regular subdivision of $\conv(\skrifta)$}\}$ form a fan in $\rr^\skrifta$.
\end{enumerate}
\end{prop}
The fan of Proposition \ref{prop:secfan}c) is called the {\em secondary fan} of $\skrifta$, and denoted $\Phi(\skrifta)$. Proposition \ref{prop:secfan}b) shows that a subdivision corresponding to a maximal cone of $\Phi(\skrifta)$ has no refinements. Hence the maximal cones correspond precisely to the primitive regular lattice triangulations of $\conv(\skrifta)$.

\section{Basic tropical geometry}\label{section:trop}
\subsection{Tropical hypersurfaces}
The purpose of this section is to recall the basics about tropical hypersurfaces and their dual subdivisions. Good references for proofs and details are \cite{RGST}, \cite{Mikh}, and \cite{Gath}. 

We work over the {\em tropical semiring} $\rr_{tr}:=(\rr,\max,+)$. Note that some authors use $\min$ instead of $\max$ in the definition of the tropical semiring. This gives a semiring isomorphic to $\rr_{tr}$. Most statements of tropical geometry are independent of this choice, but sometimes care has to be taken (cf. Lemma \ref{rem:maxmin}). 

To simplify the reading of tropical expressions, we adopt the following convention: If an expression is written in quotation marks, all arithmetic operations should be interpreted as tropical. Hence, if $x,y\in\rr$ and $k\in\nn_0$ we have for example $``x+y"=\max\{x,y\}$, $``xy"=x+y$ and $``x^k\:"=kx$. 

A {\em tropical monomial} in $n$ variables is an expression of the form $``x_1^{a_1}\dotsm x_n^{a_n}\,"$, or in vector notation, $``x^a\,"$, where $x=(x_1,\dotsc,x_n)\in\rr^n$ and $a=(a_1,\dotsc,a_n)\in\nn_0^n$. Note that $``x^a\,"=\inner{a}{x}$, the Euclidean inner product of $a$ and $x$ in $\rr^n$. A {\em tropical polynomial} is a tropical linear combination of tropical monomials, i.e.
\begin{equation}\label{trpol}
f(x)=``\sum_{a\in\skrifta} \la_a x^a\,"=\max_{a\in\skrifta}\{\la_a +\inner{a}{x}\},
\end{equation}
where $\skrifta$ is a finite subset of $\nn_0^n$, and $\la_a\in\rr$ for each $a\in\skrifta$. From the rightmost expression in \eqref{trpol} we see that as a function $\rr^n\pil\rr$, $f$ is convex and piecewise linear. The {\em tropical hypersurface} $V_{tr}(f)\sub\rr^n$ is defined to be the non-linear locus of $f\colon\rr^n\pil \rr$. Equivalently, it is the set of points $x\in\rr^n$ where the maximum in \eqref{trpol} is attained at least twice.

It is well known (see e.g. \cite{RGST} and \cite{Mikh}) that $V_{tr}(f)$ is a connected polyhedral complex of pure dimension $n-1$. As a subset of $\rr^n$, $V_{tr}(f)$ is unbounded, although some of its cells may be bounded.

We next describe the very useful duality between a tropical hypersurface $V_{tr}(f)$ and a certain lattice subdivision. With $f$ as in \eqref{trpol}, we define the {\em Newton polytope} of $f$ to be the convex hull of the exponent vectors, i.e., the lattice polytope $\conv(\skrifta)\sub\rr^n$. 
As explained in Section \ref{section:reg}, the map $a\mapsto \la_a$ induces a regular subdivision of the Newton polytope $\conv(\skrifta)$; we denote this subdivision by $\Subdiv(f)$. 

Any element $\Delta\in\Subdiv(f)$ of dimension at least 1, corresponds in a natural way to a subset $\Delta^\vee\sub V_{tr}(f)$. Namely, if the vertices of $\Delta$ are $a_1,\dotsc,a_r$, then $\Delta^\vee$ is the solution set of the equalities and inequalities
\begin{equation}\label{eq:dual}
\la_{a_1}+\inner{a_1}{x}=\dotsb=\la_{a_r}+\inner{a_r}{x}\geq \la_{b}+\inner{b}{x},\,\text{ for all $b\in\skrifta\smallsetminus \{a_1,\dotsc,a_r\}$.}
\end{equation}
That $\Delta^\vee\sub V_{tr}(f)$ follows immediately from the definition of $V_{tr}(f)$, once we observe that $r\geq 2$ (this is implied by the assumption $\dim \Delta\geq 1$). In fact, $\Delta^\vee$ is a closed cell of $V_{tr}(f)$. Moreover, we have the following theorem (see \cite{Mikh}):

\begin{theo}\label{theo:dual}
The association $\Delta\mapsto \Delta^\vee$ gives a one-to-one correspondence between the $k$-cells of $\Subdiv(f)$ and the $(n-k)$-cells of $V_{tr}(f)$, for each $k=1,2,\dotsc,n$. Furthermore, for any cells $\Delta,\Lambda\in\Subdiv(f)$ of dimensions at least 1, we have that
\begin{enumerate}
\item If $\Delta$ is a face of $\Lambda$, then $\Lambda^\vee$ is a face of $\Delta^\vee$ in $V_{tr}(f)$.
\item The affine-linear subspaces $\Aff(\Delta)$ and $\Aff(\Delta^\vee)$ are orthogonal in $\rr^n$.
\item $\Delta^\vee$ is an unbounded cell of $V_{tr}(f)$ if and only if $\Delta$ is contained in a facet of the Newton polytope of $f$.
\end{enumerate}
\end{theo}
If $C$ is a cell of $V_{tr}(f)$, we denote its corresponding cell in $\Subdiv(f)$ by $C^\vee$. The cells $C$ and $C^\vee$ are said to be {\em dual} to each other.

Theorem \ref{theo:dual} is independent of the choice of $\max$ or $\min$ as the tropical addition. However, the following lemma is not (cf. Remark \ref{rem:maxmin} below). For lack of reference, we include a proof.

\begin{lem}\label{orient}
\begin{enumerate}
\item Let $X\sub \rr^2$ be a tropical curve, and $E\in X$ a vertex. If $C$ is an edge of $X$ adjacent to $E$, then the outgoing direction vector of $C$ from $E$ is an outward normal vector of $C^\vee$ relative to $E^\vee$.
\item Let $X$ be a tropical hypersurface in $\rr^n$, where $n\geq 2$, and let $C\sub X$ be a $(n-1)$-cell adjacent to a $(n-2)$-cell $E$. If $v$ is an inward normal vector of $E$ relative to $C$, then $v$ is an outward normal vector of $C^\vee$ relative to $E^\vee$. 
\end{enumerate}
\end{lem}

\begin{proof}
a) Let $X$ be defined by the polynomial $f=``\sum_{a\in\skrifta} \la_a x^a\,"=\max_{a\in\skrifta}\{\la_a+\inner{a}{x}\}$, where $\skrifta\sub\zz^2$ is finite. Let $E$ be a vertex of $X$, and $C$ an edge of $X$ adjacent to $E$. We consider first the case where $C$ is bounded. Then $C$ has a second endpoint $F$, and $\wvec{EF}$ is a direction vector of $C$ pointing away from $E$. Dually, $C^\vee$ is the common edge of the polygons $E^\vee$ and $F^\vee$. Since we already know (by Theorem \ref{theo:dual}) that $\wvec{EF}$ is orthogonal to $C^\vee$, Lemma \ref{inout} implies that all we have to do is to show that $\inner{u}{\wvec{EF}}<0$ for some vector $u$ pointing inwards from $C^\vee$ relative to $E^\vee$. 

Let $\vv(E^\vee)=\{a_1,a_2,\dotsc,a_r\}$ be the vertices of $E^\vee$, named such that $C^\vee=\wb{a_1a_2}$. Then $u=\wvec{a_2a_3}$ points inwards from $C^\vee$ relative to $E^\vee$. We claim that $\inner{\wvec{a_2a_3}}{\wvec{EF}}<0$. To prove this, observe that the vertex $E$ satisfies the system of (in)equalities
\begin{equation}\label{Eeq}
\la_{a_1}+\inner{a_1}{E}=\la_{a_2}+\inner{a_2}{E}=\dotsb=\la_{a_r}+\inner{a_r}{E}>\la_{b}+\inner{b}{E},
\end{equation}  
for all $b\in\skrifta\smallsetminus \vv(E^\vee)$. Similarly, $F$ satisfies the relations
\begin{equation}\label{Feq}
\la_{a_1}+\inner{a_1}{F}=\la_{a_2}+\inner{a_2}{F}=\la_{c}+\inner{c}{F}>\la_{d}+\inner{d}{F},
\end{equation} 
for all $c\in \vv(F^\vee)$ and $d\in\skrifta\smallsetminus \vv(F^\vee)$. Now, in particular, \eqref{Eeq} gives $\inner{a_2}{E}-\inner{a_3}{E}=\la_{a_3}-\la_{a_2}$, while \eqref{Feq} implies (setting $d=a_3$) that $\inner{a_2}{F}-\inner{a_3}{F}>\la_{a_3}-\la_{a_2}$. Hence,
\begin{equation*}
\begin{split}
\inner{\wvec{a_2a_3}}{\wvec{EF}}&=\inner{a_3-a_2}{F-E}=\inner{a_3}{F}-\inner{a_2}{F}+\inner{a_2}{E}-\inner{a_3}{E}\\
&<\la_{a_2}- \la_{a_3}+ \la_{a_3}-\la_{a_2}=0.
\end{split}
\end{equation*}
This proves the claim, and therefore that $\wvec{EF}$ is an outward normal vector of $C^\vee$ relative to $E^\vee$.

Finally we consider the case when $C$ is unbounded. If $C$ is unbounded, then $C^\vee\sub \partial(\Delta_f)$, where $\Delta_f$ is the Newton polytope of $f$. Let $f'=``f+\la_b x^b"$, where the exponent vector $b\in\zz^2$ is chosen outside of $\Delta_f$ in such a way that $C^\vee$ is not in the boundary of $\Delta_{f'}$. If the coefficient $\la_b$ is set low enough, all elements of $\Subdiv(f)$ will remain unchanged in $\Subdiv(f')$. Furthermore, all vertices of $X$, and all direction vectors of the edges of $X$, remain unchanged in $V_{tr}(f')$. In particular, $E$ is a vertex of $V_{tr}(f')$, and its adjacent edge whose dual is $C^\vee$, has the same direction vector as $C$. Since $C^\vee$ is not in the boundary, we have reduced the problem to the bounded case above. This proves the lemma.

b) Let $\pi$ be the orthogonal projection of $\rr^n$ from $\Aff(E)$ to $\Aff(E^\vee)\simeq \rr^2$. If $C_1,\dotsc,C_r$ are the $(n-1)$-cells adjacent to $E$, then $\pi(C_1),\dotsc,\pi(C_r)$ are mapped to rays or line segments in $\Aff(E^\vee)$, with $\pi(E)$ as their common endpoint. Furthermore, if $v$ is an inward normal vector of $E$ relative to $C_i$, then $v$ is a direction vector of $\pi(C_i)$ pointing away from $\pi(E)$. The lemma now follows from the argument in a).
\end{proof}

\begin{remark}\label{rem:maxmin}
If working over the semiring $(\rr,\min,+)$ instead of $(\rr,\max,+)$, the word ``outward" in each part of Lemma \ref{orient} must be changed to ``inward".  
\end{remark}

\subsection{Tropical surfaces in $\rr^3$}
A tropical hypersurfaces in $\rr^3$ will be called simply a {\em tropical surface}. We will usually restrict our attention to those covered by the following definition:

\begin{dfn}\label{def:smooth}
Let $X=V_{tr}(f)$ be a tropical surface, and let $\delta\in\nn$. We say that the {\em degree} of $X$ is $\delta$ if the Newton polytope of $f$ is the simplex
\begin{equation*}
\Ga_\delta:=\conv(\{(0,0,0),(\delta,0,0),(0,\delta,0),(0,0,\delta)\}).
\end{equation*}
If $\Subdiv(f)$ is an elementary (unimodular) triangulation of $\Ga_\delta$, then $X$ is {\em smooth}.
\end{dfn}

\begin{remark}
We will frequently talk about a tropical surface $X$ of degree $\delta$ without referring to any defining tropical polynomial. It is then to be understood that $X=V_{tr}(f)$ for some $f$ with Newton polytope $\Ga_\delta$. In this setting, the notation $\Subdiv_X$ refers to $\Subdiv(f)$. 
\end{remark}

Let us note some immediate consequences of Definition \ref{def:smooth}. For example, since any elementary triangulation of $\Ga_\delta$ has $\delta^3$ maximal elements, $X$ must have $\delta^3$ vertices. Furthermore, any 1-cell $E\sub X$ has exactly 3 adjacent 2-cells, namely those dual to the sides of the triangle $E^\vee$. This last property makes it particularly easy to state and prove the so-called {\em balancing property}, or {\em zero-tension property} for smooth tropical surfaces. (A generalization of this holds for any tropical hypersurface. However, this involves assigning an integral {\em weight} to each maximal cell of $X$, a concept we will not need here.) 
\begin{lem}[Balancing property for smooth tropical surfaces]\label{balance}
For any 1-cell $E$ of a smooth tropical surface $X$, consider the 2-cells $C_1,C_2,C_3$ adjacent to $E$. Choosing an orientation around $E$, each $C_i$ has a unique primitive normal vector $v_i$ compatible with this orientation. Then $v_1+v_2+v_3=0$.
\end{lem}
\begin{proof}
As explained above, $C_1^\vee$, $C_2^\vee$ and $C_3^\vee$ are the sides of the triangle $E^\vee$. Theorem \ref{theo:dual} implies that $C_i^\vee$ is parallel to $v_i$ for each $i=1,2,3$. In fact, since $C_i^\vee$ is primitive, it must also have the same length as (the primitive) vector $v_i$. The vectors forming the sides of any polygon (following a given orientation), sum to zero, thus the lemma is proved.
\end{proof}

Note that when $\dim E=1$, Theorem \ref{theo:dual} guarantees that $\dim E^\vee=2$; in particular $E^\vee$ is non-degenerate. This implies that no two of the vectors $v_1,v_2,v_3$ in Lemma \ref{balance} are parallel. Thus:
\begin{lem}\label{lem:smedge}
Let $C_1,C_2,C_3$ be the adjacent 2-cells to a 1-cell of a smooth tropical surface. Then $C_1,C_2,C_3$ span different planes in $\rr^3$. 
\end{lem}

We conclude these introductory remarks on tropical surfaces with a description of some important group actions.
Let $S_4$ be the group of permutations of four elements, so that $S_4$ is the symmetry group of the simplex $\Ga_\delta$. In the obvious way this gives an action of $S_4$ on the set of subdivisions of $\Ga_\delta$.

We can also define an action of $S_4$ on the set of tropical surfaces of degree $\delta$. Let $X=V_{tr}(f)$, where $f(x_1,x_2,x_3)=``\sum_{a\in\Ga_\delta} \la_a x_1^{a_1}x_2^{a_2}x_3^{a_3}\:"$. For a given permutation $\sigma\in S_4$, we define $\sigma(X)$ as follows. First, homogenize $f$, giving a polynomial in four variables: $$f^{hom}(x_1,x_2,x_3,x_4)=``\sum_{a\in\Ga_\delta} \la_a x_1^{a_1}x_2^{a_2}x_3^{a_3}x_4^{\delta-a_1-a_2-a_3}\:".$$ Now $\sigma$ acts on $f^{hom}$ in the obvious way by permuting the variables, giving a new tropical polynomial $\sigma(f^{hom})$. Dehomogenizing again, we set $\sigma(f):=\sigma(f^{hom})(x_1,x_2,x_3,0)$. (Note that $0$ is the multiplicative identity element of $\rr_{tr}$.) Finally, we define $\sigma(X)$ to be the surface $V_{tr}(\sigma(f))$. Clearly, $\sigma(X)$ is still of degree $\delta$. The resulting action is compatible with the action of $S_4$ on the subdivisions of $\Ga_\delta$. In other words, $\Subdiv_{\sigma(X)}=\sigma(\Subdiv_X)$.

\subsection{Tropical lines in $\rr^3$}\label{section:troplines}
\begin{figure}[htbp]
\begin{center}
\input{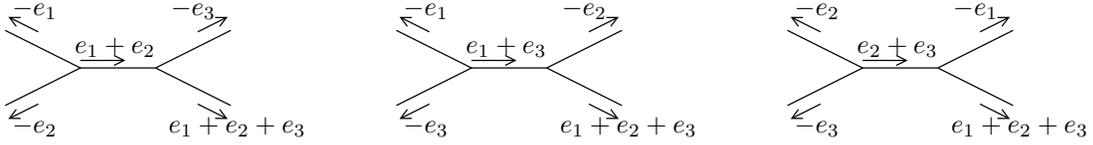}
\caption{The combinatorial types of tropical lines in $\rr^3$.}\label{spacelines} 
\end{center}
\end{figure}

Let $L$ be an unrooted tree with five edges, and six vertices, two of which are 3-valent and the rest 1-valent. We define a {\em tropical line} in $\rr^3$ to be any realization of $L$ in $\rr^3$ such that 
\begin{itemize}
\item the realization is a polyhedral complex, with four unbounded rays (the 1-valent vertices of $L$ are pushed to infinity),
\item the unbounded rays have direction vectors $-e_1$, $-e_2$, $-e_3$, $e_1+e_2+e_3$,
\item The realization is balanced at each vertex, i.e., the primitive integer vectors in the directions of all outgoing edges adjacent to a given vertex, sum to zero.
\end{itemize}
If the bounded edge has length zero, the tropical line is called {\em degenerate}. For non-degenerate tropical lines, there are three combinatorial types, shown in Figure \ref{spacelines}. From left to right we denote these combinatorial types by $(12)(34)$, $(13)(24)$ and $(14)(23)$, respectively, so that each pair of digits indicate the directions of two adjacent rays. Likewise, the combinatorial type of a degenerate tropical line is written $(1234)$.

\begin{remark}
This definition is equivalent to the more standard algebraic definition of tropical lines in $\rr^3$. See \cite[Examples 2.8 and 3.8]{RGST}. 
\end{remark}

The {\em Tropical Grassmannian}, $G_{tr}(1,3)$, is the space of all tropical lines in $\rr^3$. It is a polyhedral fan in $\rr^4$ consisting of three 4-dimensional cones, one for each combinatorial type. These cones are glued along their common lineality space of dimension 3 (corresponding to rigid translations in $\rr^3$). 

\begin{remark}
One can define tropical lines in $\rr^n$ and their Grassmannians for any $n\geq 2$. A detailed description of these spaces are given in \cite{SS}.
\end{remark}

In classical geometry, any two distinct points lie on a unique line. When we turn to tropical lines, this is true only for generic points. In fact, for special choices of points $P$ and $Q$ there are infinitely many tropical lines passing through $P$ and $Q$. The precise statement is as follows:

\begin{lem}\label{lem:twopoint}
Let $P,Q\in\rr^3$. There exist infinitely many tropical lines containing $P$ and $Q$ if and only if the coordinate vector $Q-P$ contains either a zero, or two equal coordinates. In all other cases, $P$ and $Q$ lie on a unique tropical line.
\end{lem}

An infinite collection of tropical lines in $\rr^3$, is called a {\em two-point family} if there exist two points lying on all tropical lines in the collection. Using Lemma \ref{lem:twopoint} it is not hard to see that the tropical lines of any two-point family have in fact a one-dimensional common intersection.

\section{Constructing regular elementary triangulations}\label{section:construction}
Suppose $\Delta$ is a lattice polytope contained in $\Ga_\delta$ for some $\delta\in \nn$.
We say that $\Delta$ is a {\em truncated} version of $\Ga_\delta$, if $\Delta$ results from chopping off one or several corners of $\Ga_\delta$ such that i) each chopped off piece is congruent to $\Ga_s$ for some $s<\delta$, and ii) any two chopped off pieces have disjoint interiors.   

The aim of this section is to prove that is a truncated version of $\Ga_\delta$ admits a regular, elementary triangulation (or {\em RE-triangulation} for short), then this can be extended to a RE-triangulation of $\Ga_\delta$. This fact and the lemmas building up to its proof are useful for proving existence of smooth tropical surfaces with particular properties.

We start with an easy observation, which we state in some generality for later convenience. (Recall in particular that any RE-triangulation is primitive.)
\begin{lem}\label{lem:tetr}
Suppose $\Delta\sub \rr^n$ is a $n$-dimensional lattice polytope, $F_1,F_2$ are disjoint closed faces of $\Delta$, and $\alpha_j\colon F_j\cap \zz^n\pil \rr$ is a lifting function for each $j=1,2$, such that the following conditions are fulfilled:
\begin{enumerate}[i)]
\item $\Delta=\conv(F_1\cup F_2)$,
\item $\dim(F_1)+\dim(F_2)=n-1$,
\item $\Delta\cap \zz^n=(F_1\cap \zz^n)\cup (F_2\cap \zz^n)$,
\item $\alpha_j$ induces a primitive triangulation of $F_j$, with $N_j$ maximal elements.
\end{enumerate}
Then $\alpha\colon \Delta\cap\zz^n\pil \rr$, defined by $\alpha(v):=\alpha_j(v)$ if $v\in F_j$, induces a primitive triangulation of $\Delta$ with $N_1\cdot N_2$ maximal elements, each of which is of the form $\conv(\La_1\cup \La_2)$, where $\La_i\sub F_j$ is a maximal element of the triangulation induced by $\alpha_j$.
\end{lem}
\begin{proof}
For each $j=1,2$, let $\La_j\sub F_j$ be an arbitrary maximal element of the triangulation induced by $\alpha_j$. Then $\Omega:=\conv(\La_1\cup\La_2)$ is the convex hull of $\dim(F_1)+1+\dim(F_2)+1=n+1$ lattice points, and it is a primitive simplex contained in $\Delta$. All we have to prove is that $\Omega$ is in the subdivision induced by $\alpha$. To show this, it is enough to check that 
\begin{equation}\label{eqn:obs}
\alpha(v)<\Aff_{\alpha,\Omega}(v),
\end{equation}
for any $v\in (\Delta\cap \zz^n)\smallsetminus \Omega$, where $\Aff_{\alpha,\Omega}$ is defined as the affine function extending $\alpha|_{\Omega\cap \zz^n}$ to $\Aff(\Omega)$. By condition iii), we have $v\in F_j$ for some $j=1,2$. In particular, $v\in\Aff(\La_j)$, which implies $\Aff_{\alpha,\Omega}(v)=\Aff_{\alpha_j,\La_j}(v)$. Hence \eqref{eqn:obs} is equivalent to $\alpha_j(v)<\Aff_{\alpha_j,\La_j}(v)$. But this is true since $\La_j$ is an element of the subdivision induced by $\alpha_j$. 
\end{proof}

\begin{lem}\label{lem:polyunion}
Let $\Delta_1$ and $\Delta_2$ be lattice polytopes such that $\Delta_1\cup \Delta_2$ is convex and $F:=\Delta_1\cap \Delta_2$ is a facet of both. Suppose $\ks_1$ and $\ks_2$ are regular lattice subdivisions of $\Delta_1$ and $\Delta_2$ respectively, such that $\ks_1$ and $\ks_2$ have associated lifting functions $\alpha_1$ and $\alpha_2$ which coincide on the lattice points in $F$. Then $\ks_1 \cup \ks_2$ is a regular lattice subdivision of $\Delta_1\cup \Delta_2$.
\end{lem}


\begin{proof}
Let $L(x)=0$ be the equation of the affine hyperplane spanned by $F$. For any $\la\in \rr$ consider the lifting function $\alpha_\la$ defined on the lattice points of $\Delta_1\cup \Delta_2$ by
\begin{equation*}
\alpha_\la(v):=
\begin{cases} \alpha_1(v), & \text{if $v\in \Delta_1$,}
\\
\alpha_2(v)-\la L(v), &\text{if $v \in \Delta_2$.}
\end{cases}
\end{equation*}
For $\la$ large enough, $\alpha_\la$ is concave at every point of $F$, and the induced subdivisions on $\Delta_1$ and $\Delta_2$ are $\ks_1$ and $\ks_2$ respectively.
\end{proof}

Zooming in to $\rr^3$, we now prove an auxiliary result:
\begin{lem}\label{eucl_decomp}
Let $d>e$ be natural numbers, and define the triangles $T_0,T_1 \sub \rr^3$ by 
\begin{equation*}\begin{split}
T_0&=\conv(\{(0,0,0),(d,0,0),(0,d,0)\}),\\
T_1&=\conv(\{(0,0,1),(e,0,1),(0,e,1)\}).
\end{split}\end{equation*}
Let $\kt_i$ be any RE-triangulation of $T_i$, $i=0,1$. Then there exists a RE-triangulation $\kt$ of the polytope $\Delta=\conv(T_0\cup T_1)$ such that $\kt$ extends $\kt_0$ and $\kt_1$.
\end{lem}

\begin{proof}
The strategy is as follows: We decompose $\Delta$ into three tetrahedra, find RE-triangulations of each of them, and show that these glue together to form a RE-trian\-gulation of $\Delta$. For $i=0,1$, let $\alpha_i$ be a lifting function associated to $\kt_i$, and let $\alpha\colon \Delta\cap \zz^3\pil \rr$ be defined by $\alpha(v)=\alpha_i(v)$ if $v\in T_i$. 

The decomposition of a triangular prism into three tetrahedra is well known: Let 
\begin{equation}\label{eq:eucprep}
\begin{split}
\Delta_0&=\conv(T_0\cup \{(0,0,1)\}),\\
\Delta_1&=\conv(T_1\cup \{(d,0,0)\}),\\
\Delta_2&=\Delta\smallsetminus (\Delta_0\cup \Delta_1). 
\end{split}\end{equation}
For each $i=0,1,2$, $\alpha$ restricted to $\Delta_i\cap \zz^3$ induces a primitive triangulation $\ks_i$ on $\Delta_i$. (This follows from Lemma \ref{lem:tetr}: For $\Delta_0$ and $\Delta_1$ use the decompositions indicated in \eqref{eq:eucprep}; for $\Delta_2$ take $F_1=\conv(\{(d,0,0),(0,d,0)\})$ and $F_2=\conv(\{(0,0,1),(0,e,1)\})$.) Obviously, $\ks_0$ and $\ks_1$ extend $\kt_0$ and $\kt_1$ respectively, and are elementary. Furthermore, $\ks_2$ has $de$ maximal elements, since the edges $(d,0,0)(0,d,0)$ and $(0,0,1)(0,e,1)$ are triangulated into $d$ and $e$ pieces respectively (cf. condition $iv)$ of Lemma \ref{lem:tetr}). On the other hand, $\vol(\Delta_2)=\frac16de$, so $\ks_2$ is also elementary.

Now we glue: First let $\Delta'=\Delta_0\cup \Delta_2$. Since $\ks_0$ and $\ks_2$ come from restrictions of the same lifting function, all conditions of Lemma \ref{lem:polyunion} are met, showing that $\ks_0 \cup \ks_2$ is a RE-triangulation on $\Delta'$. Also, it follows from the proof of Lemma \ref{lem:polyunion} that we can find an associated lifting function which is equal to $\alpha$ on $\Delta_2\cap \zz^3$. But then we can use Lemma \ref{lem:polyunion} again, on $\Delta=\Delta'\cup \Delta_1$. We conclude that $\ks_0\cup \ks_1\cup \ks_2$ is a RE-triangulation of $\Delta$.
\end{proof}

\begin{cor}\label{cor:extend}
Let $\Ga\sub \rr^3$ be a lattice polytope congruent to $\Ga_\delta$ for some $\delta$. Then any RE-triangulation of any of its facets can be extended to a RE-triangulation of $\Ga$.
\end{cor}
\begin{proof}
After translating and rotating, we can assume that $\Ga=\Ga_\delta$, and that the triangulated facet is the one at the bottom, i.e., $T_0$ in the above lemma. Now choose any RE-triangulation of each triangle $T_k:=\conv(\{(0,0,k),(\delta-k,0,k),(0,\delta-k,k)\})$, $k=1,\dotsc,\delta$. Lemma \ref{eucl_decomp} then implies that each layer (of height 1) $\conv(T_{k-1},T_k)$ has a RE-triangulation extending these. Finally we can glue these together one by one, as in Lemma \ref{lem:polyunion}.
\end{proof}

We now prove the main result of this section:

\begin{prop}\label{nice}
Let $\Delta$ be a truncated version of $\Ga_\delta$ for some $\delta\in \nn$. If $\kt$ is a RE-triangulation of $\Delta$, then $\kt$ can be extended to a RE-triangulation of $\Ga_{\delta}$.
\end{prop}
\begin{proof}
Each ``missing piece'' is a tetrahedron congruent to $\Ga_s$ for some integer $s<\delta$, with a RE-triangulation (induced by $\kt$) on one of its facets. Hence, by Corollary \ref{cor:extend}, each missing piece has a RE-triangulation compatible with $\kt$. By Lemma \ref{lem:polyunion}, we can glue these triangulations onto $\kt$ one by one, thus obtaining a RE-triangulation of $\Ga_\delta$. 
\end{proof}

\section{Polytopes with exits in $\Ga_\delta$}\label{section:exits}
Let $\omega_1,\omega_2,\omega_3,\omega_4$ be the vectors $-e_1,-e_2,-e_3$ and $e_1+e_2+e_3$, respectively. For any $\delta\in\nn$, and each $i=1,2,3,4$, let $F_i$ be the facet of $\Ga_\delta$ with $\omega_i$ as an outwards normal vector. For any $p\in\rr^n$, let $\ell_{p,i}$ be the unbounded ray emanating from $p$ in the direction of $\omega_i$. Hence any tropical line in $\rr^3$ with vertices $v_1$ and $v_2$, contain the rays $\ell_{v_1,i_1},\ell_{v_1,i_2},\ell_{v_2,i_3},\ell_{v_2,i_4}$ for some permutation $(i_1,i_2,i_3,i_4)$ of $(1,2,3,4)$. 
The central theme of this paper is to examine under what conditions a tropical line can be contained in a tropical surface. A simple, but crucial observation is the following:

\begin{lem}\label{lem:exitobs}
Let $C$ be a (closed) 2-cell of a tropical surface. Then,
\begin{equation*}
\ell_{p,i}\sub C\text{ for any point $p\in C$}\quad\Longleftrightarrow \quad \text{$C^\vee$ is contained in $F_i$.}
\end{equation*}
\end{lem}

Motivated by this lemma, we make the following definition:
\begin{dfn}
Let $\Delta$ be a lattice polytope contained in $\Ga_\delta$. We say that $\Delta$ has an {\em exit in the direction of $\omega_i$} if $\dim(\Delta\cap F_i)\geq 1$. If $\Delta$ has exits in the directions of $k$ of the $\omega_i$'s, we say that $\Delta$ has {\em $k$ exits}.
\end{dfn}

It is a fun task to establish how many exits different types of subpolytopes of $\Ga_\delta$ can have. We leave the proof of this lemma to the reader:
\begin{lem}\label{triangexits}
If $\delta\geq 2$, then a primitive triangle in $\Ga_\delta$ can have at most 3 exits.
\end{lem}

The case of tetrahedra with 4 exits in $\Ga_\delta$ is an interesting one, which will be important for us towards the end of the paper. Let $\kt_\delta$ be the set of all such tetrahedra. We proceed to give a classification of the elements of $\kt_\delta$, and analyze under what conditions they can be elementary. 

For any lattice tetrahedron $\Omega\sub\Ga_\delta$ we define its {\em facet distribution} $\Fac(\Omega)$ to be the unordered collection of four (possibly empty) subsets of $[4]:=\{1,2,3,4\}$ obtained in the following way: For each vertex of $\Omega$ take the set of indices $i$ of the facets $F_i$ containing that vertex. For example, if $\Omega'\sub\Ga_2$ has vertices $(0,0,0),(0,0,1),(1,1,0),(1,0,1)$, then $\Fac(\Omega')=\{\{1,2,3\},\{1,2\},\{3,4\},\{1,4\}\}$. 

A collection of four subsets of $[4]$ is called a {\em four-exit distribution} (FED) if each $i\in [4]$ appears in exactly two of the subsets. Clearly, $\Omega$ has four exits if and only if $\Fac(\Omega)$ contains a FED. (A collection $\{J_1,J_2,J_3,J_4\}$ is {\em contained} in another collection $\{J'_1,J'_2,J'_3,J'_4\}$ if (possibly after renumerating) $J_i\sub J_i'$, for all $i=1,\dotsc,4$.) For example, with $\Omega'$ as above, $\Fac(\Omega')$ contains two FEDs: $\{\{1,2,3\},\{1,2\},\{3,4\},\{4\}\}$ and $\{\{2,3\},\{1,2\},\{3,4\},\{1,4\}\}$.

Let $\ff$ be the set of all FEDs, and consider the incidence relation
\begin{equation*}
\qq\sub \kt_\delta\times \ff,\qquad \qq:=\{(\Omega,c)\:|\: \text{$c$ is contained in $\Fac(\Omega)$}\}.
\end{equation*} 
Let $\pi_1$ and $\pi_2$ be the projections from $\qq$ to $\kt_\delta$ and $\ff$ respectively. Then $\pi_1$ is obviously surjective, but not injective (for example, the last paragraph shows that $\pi_1^{-1}(\Omega')$ consists of two elements). Note that the group $S_4$ acts on $\kt_\delta$ (induced by the symmetry action on $\Ga_\delta$), on $\ff$ (in the obvious way), and on $\qq$ (letting $\sigma(\Omega,c)=(\sigma(\Omega),\sigma(c))$). Hence we can consider the quotient incidence  
\begin{equation*}
\tilde{\qq}:=\,\qq/S_4\,\sub\, \kt_\delta/S_4\times \ff/S_4,
\end{equation*}
with the projections $\tilde{\pi}_1$ and $\tilde{\pi}_2$. We claim that the image of $\tilde{\qq}$ under $\tilde{\pi}_2$ has exactly six elements, namely the equivalence classes of the following FEDs:
\begin{equation}\label{classes}
\begin{aligned}
c_1&= \{\{1,2,3\},\{1,2,4\},\{3\},\{4\}\}, \quad & c_4&= \{\{1,2,3\},\{1,2\},\{3,4\},\{4\}\}, \\
c_2&=\{\{1,2,3\},\{1,2,4\},\{3,4\},\{\:\}\},\quad& c_5&=  \{\{1,2,3\},\{1,4\},\{2,4\},\{3\}\},  \\ 
c_3&= \{\{1,2\},\{1,2\},\{3,4\},\{3,4\}\}, \quad &c_6&=  \{\{1,2\},\{1,3\},\{2,4\},\{3,4\}\}.
 \end{aligned}
\end{equation}
The proof of this claim is a matter of simple case checking: One finds that the set $\ff/S_4$ has 11 elements. In addition to the six given in \eqref{classes} there are four elements represented by FEDs of the form $\{\{1,2,3,4\},\{..\},\{..\},\{..\}\}$. These cannot be in the image of $\tilde{\pi}_2$, since no vertex lies on all four facets. Finally there is the equivalence class of $\{1,2,3\},\{1,2,3\},\{4\},\{4\}$, which corresponds to a degenerate tetrahedron. 

Now, for $\delta\in\nn$, and each $j=1\dotsc,6$, we define the following subsets of $\kt_\delta$:
 \begin{equation}
\begin{split}
\kg_\delta^j&:=\{\Omega\in \kt_\delta\:|\:  \tilde{\Omega}\in\tilde{\pi}_1(\tilde{\pi}_2^{-1}(c_j))\}\\
 \ke_\delta^j&:=\{\Omega\in \mathcal{G}_\delta^j\:|\: \text{$\Omega$ is elementary}\}.
\end{split} \end{equation}
(Here, $\tilde{\Omega}$ denotes the image of $\Omega$ in $\kt_\delta/S_4$.) Note that for a fixed $\delta$, the subsets $\kg_\delta^j$ cover $\kt_\delta$, but may overlap. For instance, our running example $\Omega'$ lies in $\kg_2^4\cap \kg_2^6$.

In the particular case $\delta=1$, we have trivially that for all $j=1,\dotsc,6$, the sets $\kg_1^j$ and $\ke_1^j$ both consist of the single tetrahedron $\Ga_1$. For higher values of $\delta$, we have the following results for the subsets $\ke_\delta^j$: 

 \begin{prop}\label{classprop}
 Let $\delta\geq 2$ be a natural number. Then
 \begin{enumerate}
 \item $\ke_\delta^1=\ke_\delta^2=\ke_\delta^3=\emptyset$.
 \item $\ke_\delta^4\cap\ke_\delta^5\neq \emptyset$.
 \item $\ke_\delta^5\smallsetminus(\ke_\delta^4\cup\ke_\delta^6)=\emptyset$.
 \item $\ke_\delta^6\smallsetminus(\ke_\delta^4\cup\ke_\delta^5)=\emptyset\Longleftrightarrow$ either $\delta=3$, or $\delta$ is even and contained in a certain sequence, starting with $2,4,6,8,14,16,18,20,26,30,56,76,\dotsc$.
 \end{enumerate}
 \end{prop}

 \begin{proof}
 a) 
 Any tetrahedron $\Omega$ in $\mathcal{G}_\delta^1$ or $\mathcal{G}_\delta^2$ contains a complete edge of $\Ga_\delta$. Such an edge is not primitive when $\delta>1$, hence $\Omega$ cannot be elementary. 
If $\Omega\in\mathcal{G}_\delta^3$, then (mod $S_4$) the vertices of $\Omega$ are of the form $(0,0,a),(0,0,b),(c,\delta-c,0),(d,\delta-d,0)$. Its volume is
 $$\Bigl|\frac16\mbox{\scriptsize 
 $\begin{vmatrix}
 0&0&a&1\\
 0&0&b&1\\
 c&\delta-c&0&1\\
 d&\delta-d&0&1
 \end{vmatrix}$\normalsize}\Bigr|
 =|\frac16\delta(a-b)(c-d)|,$$
 which is either equal to $0$ or $\geq \frac{\delta}6$. Hence $\Omega$ cannot be elementary when $\delta>1$.\\[.1cm]
 b) 
 Given any natural number $\delta$, let $\Omega$ be the convex hull of $(0,0,0),(1,0,0),(\delta-1,0,1)$ and $(0,1,\delta-1)$. Then $\Omega\in\mathcal{G}_\delta^4\cap\mathcal{G}_\delta^5$. Also, $\vol(\Omega)=\frac16$, so $\Omega$ is elementary.\\[.1cm]
 c)
 Any $\Omega\in\mathcal{G}_\delta^5$ has (modulo $S_4$) vertices with coordinates $(0,0,0)$, $(\delta-a,0,a)$, $(0,b,\delta-b)$, and $(c,d,0)$, where $a,b,c,d$ are natural numbers such that $0\leq a,b,c,d\leq \delta$ and $c+d\leq \delta$. Furthermore, if $\Omega\notin \mathcal{G}_\delta^j$ for all $j\neq 5$, then all these inequalities are strict. If $\Omega$ is elementary, we must have $\vol(\Omega)=\frac16$, which is implies that
 \begin{equation}
   \label{t123}
 6\vol(\Omega)=\Bigl|\mbox{\scriptsize
 $\begin{vmatrix}
 \delta-a&0&a\\
 0&b&\delta- b\\
 c&d&0
 \end{vmatrix}$\normalsize}\Bigr|
 =|abc+(\delta-a)(\delta-b)d|
 \end{equation}
 is equal to 1. This is impossible when $\delta\geq 2$, as shown in Lemma \ref{diofant} below.\\[.1cm]
 d) 
The vertices of $\Omega\in\mathcal{G}_\delta^6\smallsetminus(\mathcal{G}_\delta^1\cup\mathcal{G}_\delta^2\cup\mathcal{G}_\delta^3\cup\mathcal{G}_\delta^4\cup\mathcal{G}_\delta^5)$ are (modulo $S_4$) of the form $(a,0,0)$, $(0,b,0)$, $(0,c,\delta-c)$, and $(d,0,\delta-d)$, where $1\leq a,b,c,d\leq \delta-1$. We find
 $$6\vol(\Omega)=|ac(\delta-b-d)-bd(\delta-a-c)|=\colon f(\delta,a,b,c,d).$$
 When $\delta=3$, it is straightforward to check by hand that the equation $f(\delta,a,b,c,d)=1$ has no solutions in the required domain. However, if $\delta=2n+1$ for any $n\geq 2$, then $(a,b,c,d)=(n-1,n,n+1,n)$ is a solution, since $f(2n+1,n-1,n,n,n+1)=|(n-1)(n+1)-n^2|=1$.

 When $\delta$ is even we do not have any general results. A computer search shows that the equation $f(\delta,a,b,c,d)=1$ has solutions (in the allowable domain) for all $\delta$ less than 1000 except for $\delta\in\{2, 4, 6, 8, 14, 16, 18, 20, 26, 30, 56, 76\}$. It would be interesting to know whether more exceptions exist.
 \end{proof}

 \begin{lem}\label{diofant}
 The  equation $$abc+(\delta-a)(\delta-b)d=\pm 1$$
 has no integer solutions in the domain $1\leq a,b\leq \delta-1$, $\:c,d\neq 0$.
 \end{lem}
 \begin{proof}
 Keep $c,d\in \zz\smallsetminus \{0\}$ and $\delta\in \nn$ fixed, and let $\epsilon$ be either $1$ or $-1$. Then the equation $cxy+d(\delta-x)(\delta-y)=\epsilon$ describes a hyperbola $C$ intersecting the $x$-axis in $x^*=(\delta-\frac{\epsilon}{d\delta},0)$ and the $y$-axis in $y^*=(0,\delta-\frac{\epsilon}{d\delta})$. Observe that $\delta-\frac{\epsilon}{d\delta}$ is strictly bigger than $\delta-1$, and furthermore that the slope $y'(x)=\frac{d(\delta-y)-cy}{cx-d(\delta-x)}$ is positive at both $x^*$ and $y^*$. It follows that $C$ never meets the square $1\leq x,y\leq \delta-1$. This proves the lemma.
 \end{proof}

\section{Properties of tropical lines on tropical surfaces}\label{section:genprop}
From now on, unless explicitly stated otherwise, $X$ will always be a smooth tropical surface of degree $\delta$ in $\rr^3$, and $L$ a tropical line in $\rr^3$. We fix the notation $\ell_1,\dotsc,\ell_4$ for the unbounded rays of $L$ in the directions $-e_1,-e_2,-e_3$ and $e_1+e_2+e_3$, respectively, and $\ell_5$ the bounded line segment. 

Any tropical surface $X$ induces a map $c_X$ from the underlying point set of $X$ to the set of cells of $X$, mapping a point on $X$ to the minimal cell (w.r.t. inclusion) on $X$ containing it. In particular we introduce the following notion: If $v$ is a vertex of $L\sub X$, and $\dim c_X(v)=k$, we say that $v$ is a {\em $k$-vertex} of $L$ (on $X$).

An important concept for us is the possibility of a line segment on $X$ to pass from one cell to another. When $X$ is smooth, it turns out that this can only happen in one specific way, making life a lot simpler for us. We prove this after giving a precise definition:

\begin{dfn}\label{def:tres}
Let $X$ be a tropical surface (not necessarily smooth), and let $\ell\sub X$ be a ray or line segment. Let $\cc_X(\ell)$ be the set
\begin{equation*}
\cc_X(\ell):=\{c_X(p)\:|\: \text{$p\in\ell$, and $c_X(q)=c_X(p)$ for all $q\in\ell$ sufficiently close to $p$}.\}.
\end{equation*}
If $|\cc_X(\ell)|\geq 2$, then we say that $\ell$ is {\em trespassing} on $X$.
\end{dfn}
Note that $\cc_X(\ell)$ consists of the cells $C\sub X$ which satisfy $\dim(\inter(C)\cap \ell)\geq 1$. Thus Definition \ref{def:tres} corresponds well to the intuitive concept of ``passing from one cell to another''. 

\begin{lem}\label{segmentpass}
Suppose $X$ is smooth, $\ell\sub X$ a trespassing line segment, and $C,C'\sub X$ cells such that $$\cc_X(\ell)=\{C,C'\}.$$ Then $C$ and $C'$ are maximal cells of $X$ whose intersection is a vertex of $X$.
\end{lem}
\begin{proof}
Let $E=C\cap C'$, and let $v$ be a direction vector of $\ell$. Clearly, $\dim E$ is either 1 or 0. If $E$ is a 1-cell, then $C$ and $C'$ are 2-cells adjacent to $E$. But since $X$ is smooth, Lemma \ref{lem:smedge} implies that $\ell$ cannot intersect the interiors of both $C$ and $C'$, contradicting that $\cc_X(\ell)=\{C,C'\}$. 

Hence $\dim E=0$, i.e., $E$ is a vertex of $X$. Since $X$ is smooth, $E^\vee$ is a tetrahedron in $\Subdiv_X$. Now, if $\dim C=\dim C'=1$, then both $C$ and $C'$ are parallel to $v$, implying that $E^\vee$ has two parallel facets ($C^\vee$ and $C'^\vee$). This contradicts that $E^\vee$ is a tetrahedron. The case where $\dim C=1$ and $\dim C'=2$ (or vice versa) is also impossible. Here, $C^\vee$ and $C'^\vee$ would be, respectively, a facet and an edge of $E^\vee$, where $v$ is the normal vector of $C^\vee$ and $v$ also is normal to $C'^\vee$ (since $C'^\vee$ is normal to $C'$ which contains $\ell$). This would lead to $E^\vee$ being degenerate. The only possibility left is that $\dim C=\dim C'=2$, in other words that $C$ and $C'$ are both maximal. This proves the lemma.
\end{proof}

In the following, we will call a tropical line $L$ trespassing on $X$, if $L\sub X$, and at least one of the edges of $L$ is trespassing. Obviously, Lemma \ref{segmentpass} implies that:
\begin{cor}\label{trespassvert}
Any trespassing tropical line on $X$ contains a vertex of $X$.
\end{cor}
\begin{proof}
By definition, a trespassing tropical line on $X$ has a trespassing edge (either a ray or a line segment). Then we can find a line segment $\ell$ contained in this edge, such that $|\cc_X(\ell)|=2$. By Lemma \ref{segmentpass}, $\ell$ contains a vertex of $X$.
\end{proof}

\begin{lem}\label{edgedegen}
Suppose $L\sub X$ is non-degenerate, and that $L$ has a 1-vertex $v$ on $X$. Let $E=c_X(v)$. Then we have:
\begin{enumerate}
\item $E$ contains no other points of $L$.
\item The edges of the triangle $E^\vee\sub \Subdiv_X$ are orthogonal to the vectors $\omega_i$, $\omega_j$ and $\omega_i+\omega_j$ (in some order), where $\omega_i$ and $\omega_j$ are the directions of the unbounded edges of $L$ adjacent to $v$.  
\end{enumerate}
\end{lem}
\begin{proof}
a) Since $L$ is non-degenerate, $v$ has exactly three adjacent edges. Let $m_1,m_2,m_3$ be the intersections of these with a neighborhood of $v$, small enough so that each $m_i$ is contained in a closed cell of $X$. It is sufficient to prove that none of these segments are contained in $E$. Assume otherwise that $m_1\sub E$. Since $v\in\inter(E)$, the only other cells of $X$ meeting $v$ are the three (since $X$ is smooth) 2-cells adjacent to $E$. Hence $m_2\sub C$ and $m_3\sub C'$, where $C$ and $C'$ are 2-cells adjacent to $E$. We must have $C\neq C'$, otherwise $L$ cannot be balanced at $v$. But then, since $X$ is smooth, $C$ and $C'$ span different planes in $\rr^3$ (see Lemma \ref{lem:smedge}). This again contradicts the balancing property of $L$ at $v$. Indeed, balance at $v$ immediately implies that the plane spanned by $m_1$ and $m_2$ equals the plane spanned by $m_1$ and $m_3$. \\
b) Follows from a) and Lemma \ref{lem:smedge}.
\end{proof}

\begin{cor}\label{cor:deg}
Let $v_1$ and $v_2$ be the (possibly coinciding) vertices of $L\sub X$, and let $V_i=c_X(v_i)$ for $i=1,2$. Then $L$ is degenerate if and only if $V_1=V_2$.
\end{cor}
\begin{proof}
One implication is true by definition. For the other implication, suppose $V_1=V_2=:V$. If $\dim V=0$, then $L$ is clearly degenerate. If $\dim V=1$, then we must have $v_1=v_2$ (indeed, $v_1\neq v_2$ would contradict Lemma \ref{edgedegen}a)), thus $L$ is degenerate. Finally, $\dim V$ cannot be 2, as this would imply the absurdity that $V$ spans $\rr^3$. 
\end{proof}

We are now ready to prove the following proposition:
\begin{prop}\label{prop:vertex}
If $\deg X\geq 3$, then any tropical line $L\sub X$ passes through at least one vertex of $X$.
\end{prop}

\begin{proof}
Suppose $L\cap X^0=\emptyset$. By Corollary \ref{trespassvert}, $L$ must be non-trespassing. Also, $L$ cannot be degenerate. Indeed, if it were, let $v$ be its vertex. Then $c_X(v)^\vee$ would have to be a primitive triangle in $\Ga_\delta$ with four exits, contradicting Lemma \ref{triangexits}. For non-degenerate tropical lines, it is easy to rule out all cases except for one, namely when both of $L$'s vertices are 1-vertices (necessarily on different edges on $X$), as suggested to the left in Figure \ref{A11}. We can assume w.l.o.g. that the combinatorial type $L$ is $((1,2),(3,4))$. Applying Lemma \ref{edgedegen}b), it is clear that $\Subdiv_X$ contains two triangles with a common edge, with exits as shown to the right in Figure \ref{A11}. The points $A,B,C,D$ lie on $F_{14},F_{23},F_{12},F_{34}$ respectively, and the middle edge $AB$ is orthogonal to $e_1+e_2$. It follows that the points are situated as in Figure \ref{A11sub}, with coordinates of the form $A=(a,0,0)$, $B=(0,a,\delta-a)$, $C=(0,0,c)$ and $D=(d,\delta-d,0)$. Since $X$ is smooth, the triangles $ABC$ and $ABD$ must be facets of some elementary tetrahedra $ABCP$ and $ABDQ$. Setting $P=(p_1,p_2,p_3)$ and $Q=(q_1,q_2,q_3)$ we find that
\begin{equation*}
6\vol(ABCP)=\biggl|\mbox{\scriptsize
$\begin{vmatrix}
a & 0 & 0 & 1\\
0& a& \delta-a & 1\\
0&0& c& 1\\
p_1&p_2 &p_3 &1
\end{vmatrix}$\normalsize}\biggr|=
|a(ac+\delta p_2-ap_2-ap_3-c_2-cp_1)|,
\end{equation*}
implying that $a=1$, and similarly that
\begin{equation*}
6\vol(ABDQ)
=|(\delta-a)(da-\delta a+aq_2+aq_3+\delta q_1-dq_2-dq_1)|,
\end{equation*}
giving $\delta-a=1$. Hence we conclude that $\delta=2$, as claimed.
\end{proof}

\begin{figure}[tbp]
\begin{minipage}[b]{.67\linewidth}\begin{center}
\input{A11.pstex_t}
\caption{A tropical line not containing any vertices of $X$.}\label{A11} 
\end{center}\end{minipage}\hfill
\begin{minipage}[b]{.27\linewidth}\begin{center}
\input{A11sub.pstex_t}
\caption{Positions of $A,B,C,D\in \Ga_\delta$.}
\label{A11sub}
\end{center}\end{minipage}
\end{figure}

\section{Tropical lines on smooth tropical quadric surfaces}\label{section:quad}
The aim of this section is to prove a tropical analogue of the following famous theorem in classical geometry: A smooth algebraic surface of degree two has two rulings of lines. 

We begin by describing the compact maximal cells of a smooth tropical quadric. It turns out that there is always exactly one such cell:
\begin{prop}\label{onebounded}
A smooth tropical quadric surface has a unique compact 2-cell. This cell has a normal vector of the form $-e_i+e_j+e_k$, for some permutation $(i,j,k)$ of the numbers $(1,2,3)$.
\end{prop}

\begin{figure}[tbp]
\begin{center}
\input{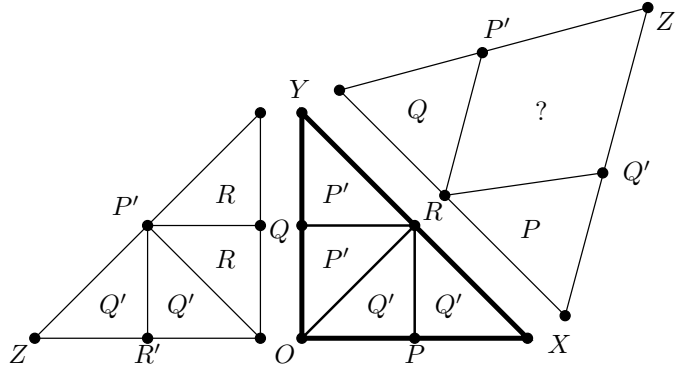}
\caption{Induced subdivisions on three facets of $\Ga_2$. A letter inside a triangle indicates the fourth point in the corresponding tetrahedron. The points $X,Y,Z,O$ are $(2,0,0),(0,2,0),(0,0,2),(0,0,0)$ respectively.}\label{sweet} 
\end{center}
\end{figure}

\begin{proof}
Let $X$ be the smooth quadric. A compact 2-cell of $X$ corresponds to a 1-cell in $\Subdiv_X$ in the interior of the Newton polytope $\Ga_2$. Such 1-cells will in the following be called {\em diagonals}.

\begin{figure}[tbp]
\begin{minipage}[b]{.43\linewidth}\begin{center}
\input{Gamma2.pstex_t}
\caption{The lattice points in $\Ga_2$.}\label{Gamma2} 
\end{center}\end{minipage}\hfill
\begin{minipage}[b]{.55\linewidth}\begin{center}
\input{facetsubs.pstex_t}
\caption{The two unique elementary triangulations of a lattice triangle with side length 2.}\label{facetsubs} 
\end{center}\end{minipage}\end{figure}

The only possible diagonals in $\Ga_2$ are the line segments (see Figure \ref{Gamma2}) 
\begin{equation}\label{diags}
PP'=\wb{(1,0,0),(0,1,1)},\quad QQ'=\wb{(1,0,1),(0,1,0)}\quad\text{and}\quad RR'=\wb{(0,0,1),(1,1,0)}.
\end{equation}
Note that all these intersect in $(\frac12,\frac12,\frac12)\notin\zz^3$, so at most one of them can be in $\Subdiv_X$. This proves uniqueness. To complete the proof we must show that $\Subdiv_X$ contains at least one diagonal. (The final statement in the proposition follows trivially from the direction vectors of the diagonals in \eqref{diags}.)

Since $X$ is smooth, $\Subdiv_X$ is an elementary triangulation of $\Ga_2$. In particular, the induced subdivisions of the four facets of $\Ga_2$ are also elementary triangulations. Up to symmetry, there are only two possibilities for these triangulations, shown as $I$ and $II$ in Figure \ref{facetsubs}. Suppose the triangulation of the bottom facet is of type $I$. Then, in particular, it contains the triangle $\triangle PQR$ as an element. Let $T\in\Subdiv_X$ be the (unique) elementary tetrahedron having this triangle as a facet. For $T$ to have volume $\frac16$, its height must be 1, so the fourth vertex is either $P'$, $Q'$ or $R'$. In either case, $T$ contains one of the diagonals \eqref{diags} as an edge.

The same argument can be used on the three other facets of $\Ga_2$, so we are left with the case where all the subdivisions induced on the facets are of type II (cf. Figure \ref{facetsubs}). Suppose this is the case, and that $\Subdiv_X$ contains no diagonals. We will show that this leads to a contradiction. 

Figure \ref{sweet} shows three of the facets of $\Ga_2$ folded out. Starting from the bottom facet $OXY$ (drawn in bold lines in Figure \ref{sweet}), we can assume (after a rotation if necessary) that its induced subdivision is as in Figure \ref{sweet}. Now, since $\Subdiv_X$ contains neither $PP'$, $QQ'$ nor $RR'$, the tetrahedron containing $OPR$ as a facet, must have $Q'$ as its fourth vertex. Similarly, the other three tetrahedra on the bottom of $\Subdiv_X$ are uniquely determined. This in turn determines the subdivision of the facet $OYZ$, and the corresponding closest tetrahedra (see Figure \ref{sweet}). In particular, it follows that $P'Q'\in\Subdiv_X$. But turning to the facet $XYZ$, we see that this is impossible. Indeed, we already know that $P'R$ and $Q'R$ are in $\Subdiv_X$. Together with $P'Q'$, this implies that the induced subdivision of $XYZ$ is of type I, violating the assumption.
\end{proof}

Let $\wb{X}$ denote the compact 2-cell of $X$ found in Proposition \ref{onebounded}. 
Our main result about tropical lines on tropical quadrics is the following:

\begin{theo}\label{theo:quad}
For each point $p\in\wb{X}$ there exist two distinct tropical lines on $X$ passing through $p$. 
\end{theo}

\begin{proof}
We can assume (using if necessary the action of $S_4$) that $\wb{X}$ has a normal vector $-e_1+e_2+e_3$, i.e., that the edge in $\Subdiv_X$ corresponding to $\wb{X}$ is $PP'$ (see Figure \ref{Gamma2}). Let $p$ be any point on $\wb{X}$, and consider the line given by $p+t(e_1+e_2)$, $t\in\rr$. Let $L^-$ and $L^+$ be the rays where $t\leq0$ and $t\geq 0$ respectively, and let $p^-, p^+$ be the points on the boundary of $\wb{X}$ where $L^-$ and $L^+$ leave $\wb{X}$. We will show that the tropical line $L_p$ with vertices $p^-$ and $p^+$, lie on $X$. 

Let $E^-:=c_X(p^-)$ and $E^+:=c_X(p^+)$. If $E^-$ (resp. $E^+$) is a vertex, redefine it to be any adjacent edge (of $\wb{X}$) not parallel to $v$. To prove that $L_p\sub X$, it is enough (by Lemma \ref{lem:exitobs}) to show that the triangle $(E^-)^\vee\in \Subdiv_X$ has exits in the directions $\omega_1,\omega_2$, and that $(E^+)^\vee$ has exits in the directions $\omega_3,\omega_4$.

The boundary of $\wb{X}$ is made up precisely by the 1-cells of $X$ whose dual triangles in $\Subdiv_X$ has $PP'$ as one edge. In particular there are lattice points $A,B\in \Ga_2$ such that $(E^-)^\vee=\triangle APP'$ and $(E^+)^\vee=\triangle BPP'$. We claim that 
\begin{equation}\label{aonf}
\text{$A$ and $B$ lies on the edges $F_{12}$ and $F_{34}$ respectively.}
\end{equation}
If this claim is true, it follows immediately that the triangles $\triangle APP'$ and $\triangle BPP'$ have the required exits, and therefore that $L_p\sub X$. To prove the claim, we utilize Lemma \ref{lem:linalg} below. By the construction of $E^-$, it is clear that the vector $e_1+e_2$ points inwards from $E^-$ into $\wb{X}$. The lemma then implies that $\inner{e_1+e_2}{u}<0$ for all vectors $u$ pointing inwards from $PP'$ into $\triangle APP'$. In particular, choosing $u$ as the vector from $P$ to $A=(a_1,a_2,a_3)$, this gives $a_1+a_2<1$. The only lattice points in $\Ga_2$ satisfying this are those on $F_{12}$, so $A\in F_{12}$. That $B\in F_{34}$ follows similarly. This proves the claim, and we conclude that $L_p\sub X$.

\begin{figure}[tbp]
  \centering
  \epsfig{file=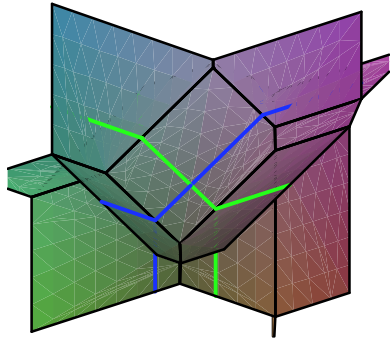,height=5cm}
  \caption{A smooth tropical quadric surface $X$, with two tropical lines passing through a point in $\wb{X}$.}\label{twolinesfig}
\end{figure}

Next, consider the affine line $p+t(e_1+e_3)$, $t\in\rr$. The points where this line leaves $\wb{X}$ are again the vertices of a tropical line, $L_p'$, which we claim is contained in $X$. Indeed, this follows after swapping the coordinates $e_2$ and $e_3$ (i.e., letting the transposition $\sigma=(23)\in S_4$ act on $X$), and repeating the above proof word for word. Figure \ref{twolinesfig} shows $L_p$ and $L_p'$ in a typical situation.
\end{proof}

\begin{figure}[tbp]
\begin{center}
\input{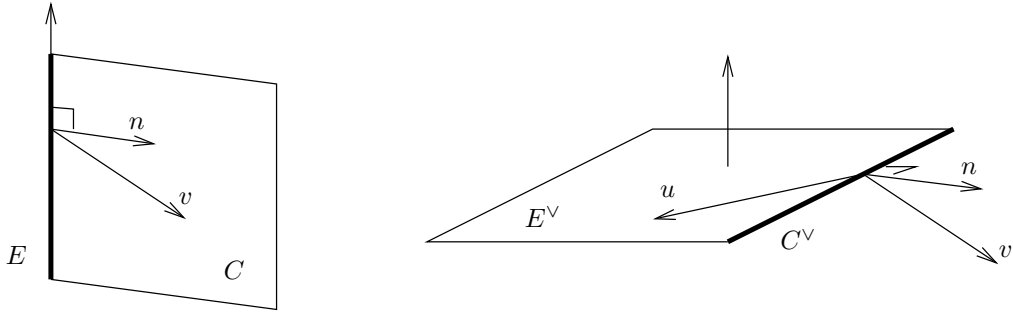}
\caption{Illustration of Lemma \ref{lem:linalg}.}\label{duallinalg} 
\end{center}
\end{figure}

\begin{lem}\label{lem:linalg}
Let $E$ be an edge of a 2-cell $C$ on a tropical surface. For any vector $v$ pointing inwards from $E$ into $C$, and any vector $u$ pointing inwards from $C^\vee$ into $E^\vee$, we have $$\inner{v}{u}<0.$$ 
\end{lem}
\begin{proof}
Let $n$ be the unit inwards normal vector of $E$ relative to $C$. By Lemma \ref{orient}, $n$ is an outwards normal vector of $C^\vee$ relative to $E^\vee$. In particular, we have $\inner{v}{n}>0$ and $\inner{u}{n}<0$. (See Figure \ref{duallinalg}.)

For $v=n$, the lemma is clearly true, so assume $v\neq n$. The vector product $v\times n$ is then a normal vector of $C$, and therefore a direction vector of $C^\vee$. Hence $u\times (v\times n)$ is a normal vector of $E^\vee$, i.e., it is a direction vector of $E$. But since $n$ is a normal vector of $E$, this implies that $\inner{u\times (v\times n)}{n}=0$. Expanding this, using the familiar formula $a\times (b\times c)=\inner{a}{c}b-\inner{a}{b}c$, we find that
\begin{equation*}
\inner{u}{n}\inner{v}{n}=\inner{u}{v}\inner{n}{n}=\inner{u}{v}.
\end{equation*}
(In the last step we used that $|n|=1$.) The lemma follows from this, since $\inner{u}{n}<0$ and $\inner{v}{n}>0$.
\end{proof}

\section{Two-point families on $X$}\label{twopoint}
To any $L\sub X$, with edges $\ell_1,\dotsc,\ell_5$, we can associate a set of data, $\dd_X(L)=\{V_1,V_2,\cc_1,\cc_2,\dotsc,\cc_5,\kappa\}$, where,
\begin{itemize}
\item $V_i=c_X(v_i)$, where $v_1,v_2$ are the (possibly coinciding) vertices of $L$.
\item $\cc_i$ is the set $\cc_X(\ell_i)$ (cf. Definition \ref{def:tres}).
\item $\kappa$ is the combinatorial type of $L$.
\end{itemize}
Recall in particular that $\ell_i$ is trespassing on $X$ if and only if $|\cc_i|\geq 2$.

\begin{figure}[htbp]
\begin{center}
\input{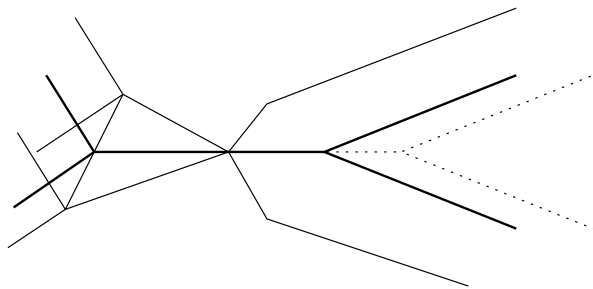}
\caption{A two-point family of tropical lines on a tropical surface.}\label{loose} 
\end{center}
\end{figure}

One might wonder if different tropical lines on $X$ can have the same set of data. It is not hard to imagine an example giving an affirmative answer, e.g. as in Figure \ref{loose}. In this Figure one of the vertices of the tropical line can be moved along the middle segment, creating infinitely many tropical lines with the same set of data. Clearly, the collection of all these tropical lines is a two-point family. As we will show in the remainder of this section, this is not a coincidence.

By a perturbation of a point $p\in\rr^3$ we mean a continuous map $\mu\colon [0,1)\pil \rr^3$, possibly constant, such that $\mu(0)=p$.
\begin{dfn}\label{def:pert}
A tropical line $L\sub X$ can be {\em perturbed on $X$} if there exist perturbations $\mu_1$ and $\mu_2$ - not both constant - of the vertices of $L$ such that for all $t\in [0,1)$, $\mu_1(t)$ and $\mu_2(t)$ are the vertices of a tropical line $L_t\sub X$.
In this case, we call the map $[0,1)\pil G_{tr}(1,3)$ given by $t\mapsto L_t$ a {\em perturbation} of $L$ on $X$.
\end{dfn}

If $L$ is degenerate, we think of $L$ as having two coinciding vertices. Thus Definition \ref{def:pert} allows perturbations of $L$ where the vertices are separated, creating non-degenerate tropical lines. 

By a {\em two-point family of tropical lines on $X$}, or simply a two-point family on $X$, we mean a two-point family of tropical lines, all of which are contained in $X$. A two-point family on $X$ is {\em maximal} (on $X$) if it not contained in any strictly larger two-point family on $X$.  A tropical line on $X$ is {\em isolated} if it does not belong to any two-point family on $X$. 

Special perturbations, as the one in Figure \ref{loose}, give rise to two-point families on $X$. We state a straightforward generalization of this example in the following lemma, for later reference. Note that if $\mu$ is a perturbation of $L$ on $X$, we say that the vertex $v_i$ is {\em perturbed along an edge} of $L$, if $\im(\mu_i)\sub\Aff(\ell)$ for some edge $\ell\sub L$ (cf. the notation in Definition \ref{def:pert}). 

\begin{lem}\label{pertfam}
If a non-degenerate $L\sub X$ has a perturbation on $X$ where at least one of the vertices is perturbed along an edge of $L$, then $L$ belongs to a two-point family on $X$. 
\end{lem}

\begin{prop}\label{prop:UD}
Let $L$ be a tropical line on a smooth tropical surface $X$, where $\deg X\geq 3$. If $L$ is isolated, then $L$ is uniquely determined by $\dd_X(L)$.
\end{prop}

\begin{proof}
Let $\dd=\dd_X(L)=\{V_1,V_2,\cc_1,\cc_2,\dotsc,\cc_5,\kappa\}$ be a given set of data. We will identify all situations where $L$ is not uniquely determined by $\dd$, and show that Lemma \ref{pertfam} applies in each of these cases.

We first consider the case where $\kappa\neq(1234)$, meaning that $L$ is non-degenerate. The following observations will be used frequently:
\begin{enumerate}[A)]
\item $L$ is determined by (the positions of) its two vertices.
\item The direction vector of the bounded segment $\ell_5$ is determined by $\kappa$.
\item If $|\cc_j|\geq 2$, then $\Aff(\ell_j)$ is determined by the elements of $C_j$ (and the index $j$).
\item If $\dim V_i=1$, and $\Aff(\ell_j)$ is known for any edge $\ell_j$ adjacent to $v_i$, then $v_i$ is determined.
\end{enumerate}
Of these, A) and B) are clear, C) is a consequence of Lemma \ref{segmentpass}, and D) follows from Lemma \ref{edgedegen}a).

Now, assume that $V_1$ and $V_2$ are ordered so that $\dim V_1\leq \dim V_2$. Under this assumption, we examine the uniqueness of $L$ for different sets of data, according to the pair $(\dim V_1,\dim V_2)$:

$\bullet\;(\dim V_1,\dim V_2)=(0,0)$: Obviously, by A), $L$ is determined.

$\bullet\;(\dim V_1,\dim V_2)=(0,1)$: In this case $\Aff(\ell_5)$ is determined by $V_1$ and $\kappa$ (cf. B)). Hence $v_2$ is determined (by D)). Since $v_1=V_1$, it follows that $L$ is determined.

$\bullet\;(\dim V_1,\dim V_2)=(0,2)$: Again, $v_1$ and $\Aff(\ell_5)$ are determined by $V_1$ and $\kappa$. Write $\kappa=((a,b),(c,d))$, and consider first the case where either $|\cc_c|\geq2$ or $|\cc_d|\geq2$. We can assume the former. Then $\Aff(\ell_{c})$ is determined, which again determines $v_2=\Aff(\ell_5)\cap \Aff(\ell_c)$. Thus, in this case $L$ is determined.

Otherwise, we have $\cc_c=\cc_d=V_2$. In this situation $L$ is not uniquely determined by $\dd$, as $v_2$ can be perturbed to anywhere in the intersection of $\Aff(\ell_5)$ and $V_2$ without changing $\dd$.

$\bullet\;(\dim V_1,\dim V_2)=(1,1)$: Observe first that we must have $|\cc_i|\geq 2$ for some $i$. (Otherwise $L$ is not trespassing, and since none of its vertices are vertices of $X$, this would contradict Proposition \ref{prop:vertex}.) 
Hence $\Aff(\ell_i)$ is determined for some $i$. If $i=5$, then (by D)) both $v_1$ and $v_2$ are determined by this. If $i\neq 5$, then in the first place only the endpoint of $\ell_i$ is determined. But this together with $\kappa$ determines $\Aff(\ell_5)$, and thus both vertices. Hence, in any case, $L$ is determined.

$\bullet\;(\dim V_1,\dim V_2)=(1,2)$: Let $\kappa=((a,b),(c,d))$. We consider five cases: 

i) $|\cc_j|\geq 2$ for both $j=c,d$. Then $\Aff(\ell_c)$ and $\Aff(\ell_d)$ are determined, and therefore also $v_2=\Aff(\ell_c)\cap \Aff(\ell_d)$. This and $\kappa$ determines $\Aff(\ell_5)$, which in turn (by D)) determines $v_1$. Hence $L$ is determined. 

ii) $|\cc_j|\geq 2$ for exactly one index $j\in \{c,d\}$ (assume $d$), and also for at least one index $j\in\{a,b,5\}$. This last condition determines $\Aff(\ell_5)$, either directly (if $j=5$) or via $v_1$ and $\kappa$. Thus $v_2=\Aff(\ell_d)\cap \Aff(\ell_5)$ is determined, and therefore $L$ as well. 

iii) $|\cc_j|\geq 2$ for exactly one index $j\in \{c,d\}$ (assume $d$), and for no other indices $j$. In this case $v_2$ can be perturbed along $\ell_d$ without changing $\dd$, so $L$ is not determined by $\dd$. (The perturbation of $v_1$ (along $V_1$) will be determined by the perturbation of $v_2$.)

iv) $|\cc_j|\geq 2$ for no $j\in \{c,d\}$, but at least one $j\in\{a,b,5\}$. As in ii) above, the last condition determines $\Aff(\ell_5)$ and therefore $v_1$. The vertex $v_2$ can be perturbed along $\ell_5$, so $L$ is not determined.

v) $|\cc_j|=1$ for all $j\in \{1,2,3,4,5\}$. This is not possible when $\deg X\geq 3$. In fact, it follows from Lemma \ref{triangexits} that $\deg X=1$. Indeed, since no edge of $L$ is trespassing, the triangle $V_1^\vee$ must have four exits in $\Ga_{\deg X}$.

$\bullet\;(\dim V_1,\dim V_2)=(2,2)$: Note first that $V_1\neq V_2$, since $L$ spans $\rr^3$. Hence $|\cc_5|\geq 2$, determining $\Aff(\ell_5)$. Now, for both $i=1,2$ we have: If any adjacent unbounded edge of $v_i$ is trespassing, then $v_i$ is determined. If not, $v_i$ can be perturbed along $\ell_5$ keeping $\dd$ unchanged.

Going through the above list, we see that in each case where $L$ is not uniquely determined by $\dd$, $L$ has a perturbation where a vertex is perturbed along an edge of $X$. Hence, by Lemma \ref{pertfam}, $L$ belongs to a two-point family on $X$. 

Finally, suppose $\kappa=(1234)$, so $L$ is degenerate. We show that in this case, $L$ is determined by $\dd$. Corollary \ref{cor:deg} (and its proof) tells us that $V_1=V_2:=V$ where $\dim V$ is either 0 or 1. In the first case, $L$ is obviously uniquely determined. If $\dim V=1$ then $|\cc_j|\geq 2$ for some $j\in\{1,2,3,4\}$, otherwise $L$ would contain no vertex of $X$, contradicting Proposition \ref{prop:vertex}. Hence $\Aff(\ell_j)$ is determined. We claim that $V_1\not\sub\Aff(\ell_j)$. Note that this would determine $v_1=v_2=\Aff(\ell_j)\cap V_1$, and therefore also $L$. To prove the claim, note that if $V_1\sub\Aff(\ell_j)$, then $V_1\in \cc_j$. This is impossible, since any element of $\cc_j$ must be of dimension 2 (cf. Lemma \ref{segmentpass}). This concludes the proof of the proposition.
\end{proof}

\section{Tropical lines on higher degree tropical surfaces}\label{section:higher}
In this section we present our main results about tropical lines on smooth tropical surfaces of degree greater than two. The proofs rest heavily on what we have done so far. The  first is indeed a corollary of Proposition \ref{prop:UD}:

\begin{cor}
Let $X$ be a smooth tropical surface where $\deg X\geq 3$. Then $X$ contains at most finitely many isolated tropical lines. Furthermore, $X$ contains at most finitely many maximal two-point families.
\end{cor}
\begin{proof}
The first statement is immediate from Proposition \ref{prop:UD}, since there are only finitely many possible sets of data $\dd_X(L)$. For the last statement, observe that any two-point family contains a non-degenerate tropical line. Going through the proof of Proposition \ref{prop:UD}, we see that if $\dd$ is the data set of is a non-degenerate tropical line, then there can be at most one maximal two-point family containing tropical lines with data set $\dd$. Hence there are at most finitely many maximal two-point families on $X$.
\end{proof}

The next theorem show that two-point families exist on smooth tropical surfaces of any degree.
\begin{theo}\label{fullcone}
For any integer $\delta$, there exists a full dimensional cone in $\Phi(\Ga_\delta)$ in which each point corresponds to a smooth tropical surface containing a two-point family of tropical lines. In particular, there exist smooth tropical surfaces of degree $\delta$ with infinitely many tropical lines.
\end{theo}

\begin{figure}[tbp]
\begin{minipage}[b]{.48\linewidth}
  \centering
  \epsfig{file=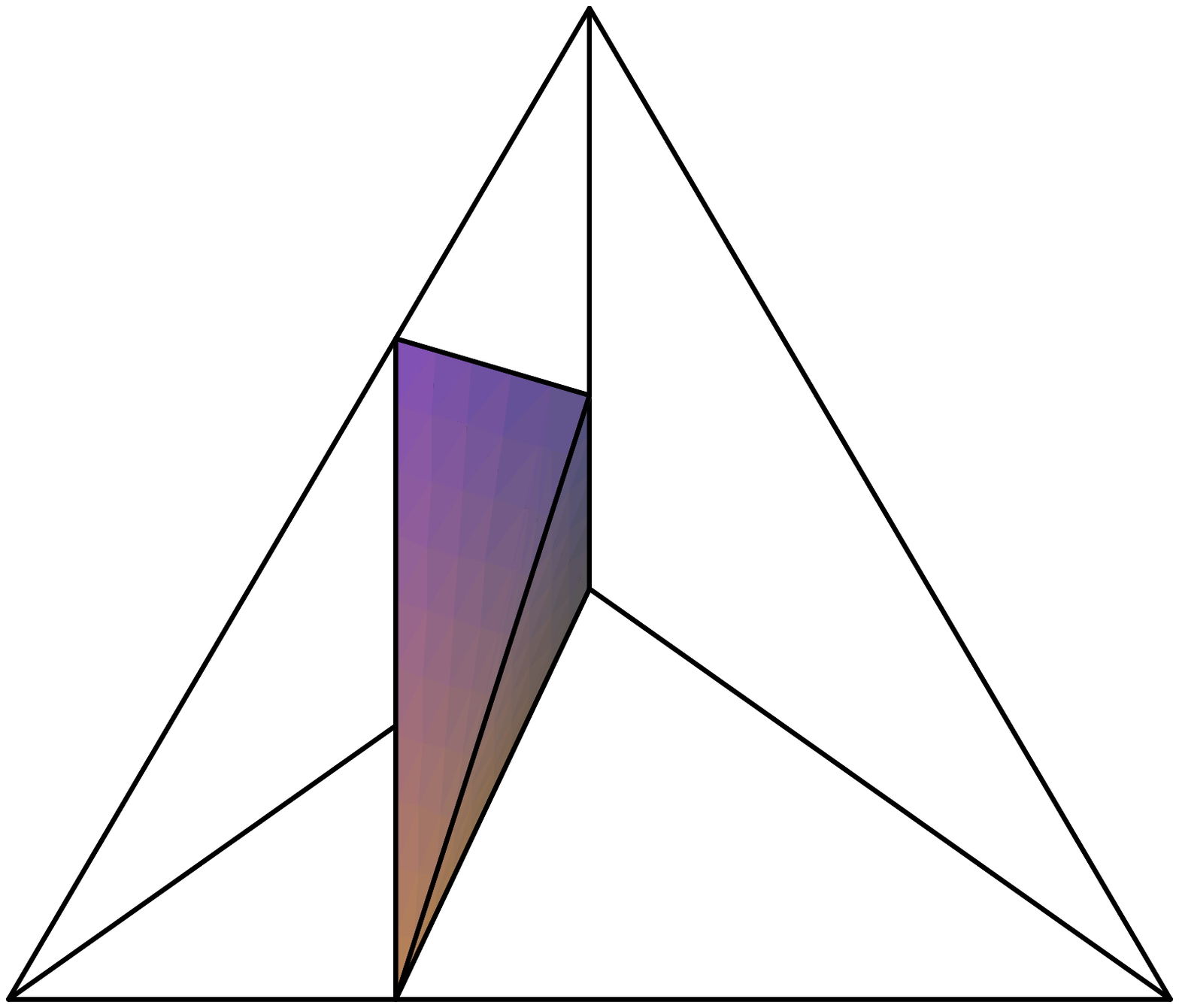,height=4.2cm}
  \caption{A tetrahedron with four exits in $\Ga_3$.}\label{classIV}
\end{minipage}\hfill
\begin{minipage}[b]{.48\linewidth}
\centering
  \epsfig{file=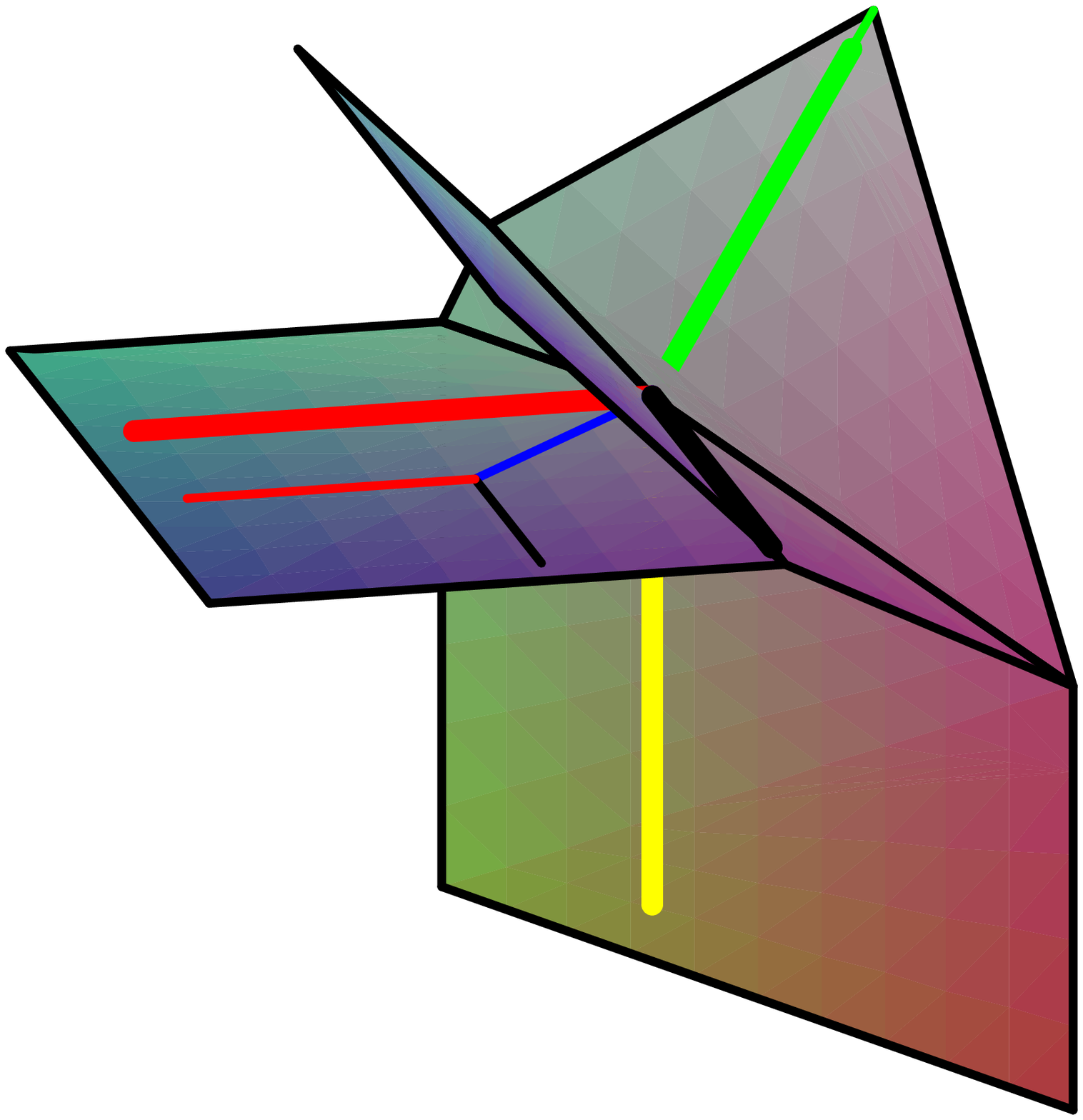,height=4cm}
  \caption{The degenerated tropical line corresponding to the tetrahedron in Figure \ref{classIV} belongs to a two-point family}\label{loosenable}
\end{minipage}
\end{figure}

\begin{proof}
Let $\delta$ be an arbitrary, fixed integer. Consider the lattice tetrahedron $\Omega\sub\rr^3$ defined by 
\begin{equation}\label{eq:omegatetr}
\Omega_\delta:=\conv(\{(0,0,0),(0,0,1),(\delta-1,1,0),(1,0,\delta-1)\}).
\end{equation}
It is easy to see that $\Omega_\delta$ has four exits in $\Ga_\delta$ (see Figure \ref{classIV}).

Assume for the moment that there exists a smooth tropical surface $X$ of degree $\delta$ such that $\Subdiv_X$ contains $\Omega_\delta$. Then Lemma \ref{lem:exitobs} implies the vertex $v:=\Omega_\delta^\vee\in X$ is the center of degenerate tropical line $L\sub X$. We claim that $L$ belongs to a two-point family on $X$. Indeed, this also follows from Lemma \ref{lem:exitobs}: Let $C\sub X$ be the cell dual to the line segment in $\Subdiv_X$ with vertices $(0,0,0)$ and $(0,0,1)$. Then for any point $p(t)=v+t(-e_1-e_2)$, where $t>0$, the line segment with endpoints $v$ and $p(t)$ is contained in $C$. Let $L_t$ be the tropical line with vertices $v$ and $p(t)$. Lemma \ref{lem:exitobs} guarantees that the rays starting in $p(t)$ in the directions $-e_1$ and $-e_2$ are contained in $C$. Hence $L_t\sub X$. Clearly, the lines $L_t$ form a two-point family on $X$, thus the claim is true. (See Figure \ref{loosenable}.)

It remains to prove the existence of a RE-triangulation of $\Ga_\delta$ containing $\Omega_\delta$. Using the techniques in Section \ref{section:construction}, it is not hard to construct such a triangulation explicitly. For example, consider the polytope 
$$\Delta=\conv(\{(0,0,0),(\delta,0,0),(\delta-1,1,0),(0,1,0),(0,1,\delta-1),(0,0,\delta)\}).$$ Then $\Delta$ is a truncated version of $\Ga_\delta$, so by Proposition \ref{nice} it is enough to construct a RE-triangulation of $\Delta$ which contains $\Omega_\delta$. Write $\Delta=\Omega_\delta\cup \Delta_1\cup \Delta_2\cup \Delta_3\cup \Delta_4$, where
\begin{align*}
  \Delta_1&=\conv(\{(0,0,0),(\delta,0,0),(\delta-1,1,0),(1,0,\delta-1)\})\\
  \Delta_2&=\conv(\{(0,0,1),(\delta-1,1,0),(1,0,\delta-1),(0,0,\delta)\})\\
  \Delta_3&=\conv(\{(0,0,0),(\delta-1,1,0),(0,1,0),(0,0,\delta)\})\\
  \Delta_4&=\conv(\{(\delta-1,1,0),(0,1,0),(0,1,\delta-1),(0,0,\delta)\}).
\end{align*}
Repeated use of Lemma \ref{lem:tetr} gives a RE-triangulation of each of these (for $\Delta_1$ and $\Delta_4$ choose any RE-triangulation of the facets $\conv(\{(0,0,0),(\delta,0,0),(1,0,\delta-1)\})$ and $\conv(\{(\delta-1,1,0),(0,1,0),(0,1,\delta-1)\})$ respectively). 
Finally, it is easy to check that these triangulations patch together to a RE-triangulation of $\Delta$, using Lemma \ref{lem:polyunion}.
\end{proof}

\begin{ex}
Define the tropical polynomial $g_3$ by
\begin{multline*}
g_3(x,y,z)=``-22x^3+16x^2y-10x^2z+0xy^2+0xz^2+8xyz-23y^3-12y^2z\\
-5yz^2+0z^3-14x^2+14xy-3xz-6y^2+4yz+0z^2-8x+6y-z-3".
\end{multline*}
The subdivision $\Subdiv(g_3)$ (shown in Figure \ref{fig:subg3}) is a RE-triangulation of $\Ga_3$ containing the tetrahedron $\Omega_3$ (defined in \eqref{eq:omegatetr}). Hence $V_{tr}(g_3)$ is a smooth tropical cubic surface with a two-point family of tropical lines, all of which have $\Omega_3^\vee=(1,-21,-2))$ as a vertex. The polynomial $g_3$ was constructed by first building the RE-triangulation (following the suggestions in the proof of Theorem \ref{fullcone}, making appropriate choices where needed), and then calculating an interior point in the secondary cone of this subdivision. The latter part was done using the Maple package Convex (\cite{maple, convex}).
\begin{figure}[tbp]
  \centering
  \epsfig{file=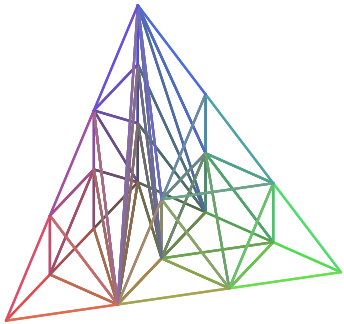,height=4.5cm}
  \caption{The RE-triangulation induced on $\Ga_3$ by $g_3$.}\label{fig:subg3}
\end{figure}

Similarly, the tropical polynomial $g_4$ below gives a smooth tropical surface of degree four containing a two-point family of tropical lines:
\begin{multline*}
g_4(x,y,z)=
``-12x^4+72x^3y-x^3z-4x^2y^2+41x^2yz+7x^2z^2-91xy^3\\
-39xy^2z+2xyz^2+12xz^3-189y^4-133y^3z-85y^2z^2-45yz^3-6z^4\\-5x^3+56x^2y+5x^2z-24xy^2+24xyz+11xz^2-118y^3-63y^2z\\-19yz^2-3z^3-x^2+32xy+7xz-55y^2-4yz-z^2+0x+0y+0z+0".
\end{multline*}
\end{ex}

In light of the above theorem, one might ask whether there exist tropical surfaces of high degree containing an {\em isolated} degenerate tropical line $L$. If we add the requirement that $L$ is non-trespassing on $X$, we can give the following partial answer:
\begin{prop}\label{prop:deglist}
Let $\delta\in \nn$. There exists a smooth tropical surface $X$ of degree $\delta$ containing an isolated, non-trespassing, degenerate tropical line, if and only if $\delta$ is
 \begin{itemize}
 \item an odd number greater than 3, or
 \item an even number except 2, 4, 6, 8, 14, 16, 18, 20, 26, 30, 56, 76,... 
 \end{itemize} 
 \end{prop}
 \begin{proof}
 We know that the vertex of such a line must be a vertex of $X$, corresponding to an elementary tetrahedron $\Omega\in \Subdiv_X$ with four exits. Furthermore, no edge of $\Omega$ can have more than one exit. Indeed, an edge with exits $\omega_i$ and $\omega_j$ will be orthogonal to the vector $\omega_i+\omega_j$, implying (as in the proof of Theorem \ref{fullcone}) that $L$ belongs to a two-point family. 

From the classification in \eqref{classes} of tetrahedra with four exits in $\Ga_\delta$, we observe the following: A tetrahedron with four exits, in which no edge has more than one exit, must belong either exclusively to the subset $\mathcal{G}_\delta^5$, or exclusively to the subset $\mathcal{G}_\delta^6$. The result then follows from Proposition \ref{classprop}c) and d). As we remarked in that proposition, we do not know how (or if) the list of even degree exceptions continues.
 \end{proof} 

Both Theorem \ref{fullcone} and Proposition \ref{prop:deglist} show that there exist plenty of tropical surfaces of arbitrarily high degree containing tropical lines. It is natural to wonder whether there also exist smooth tropical surfaces containing {\em no} tropical lines, isolated or not. This is indeed true in all degrees greater than three, as we prove in \cite{Vig_Class}. In that paper we present a classification of tropical lines on general smooth tropical surfaces, and propose a method for counting the isolated tropical lines on such surfaces. 

\vspace{.3cm}

\noindent
{\it Acknowledgements}.
I would like to thank my supervisor Kristan Ranestad for many interesting discussions about tropical geometry, and for numerous useful suggestions and improvements during the writing of this paper.

\bibliographystyle{plain}
\bibliography{bib_LoS}
\end{document}